\documentclass[iicol, sn-apa]{sn-jnl}

\usepackage{xurl}
\usepackage{verbatim}

\usepackage{optidef}
\usepackage{bm}
\usepackage[font=small]{caption}
\usepackage{subcaption}

\usepackage{graphicx}%
\usepackage{multirow}%
\usepackage{amsmath,amssymb,amsfonts}%
\usepackage{amsthm}%
\usepackage[title]{appendix}%
\usepackage{xcolor}%
\usepackage{textcomp}%
\usepackage{manyfoot}%
\usepackage{booktabs}%
\usepackage{algorithm}%
\usepackage{algorithmicx}%
\usepackage{algpseudocode}%
\usepackage{listings}%

\renewcommand{\APACrefURL}[1]{\href{https://#1}{\url{#1}}}

\usepackage{tikz}
\usetikzlibrary{shapes,positioning}

\theoremstyle{thmstyleone}%

\theoremstyle{thmstyletwo}%

\theoremstyle{thmstylethree}%

\DeclareMathOperator*{\argmin}{arg\,min}

\begin{document}

\title{\LARGE \bf
Fast Near Time-Optimal Motion Planning for Holonomic Vehicles in Structured Environments
}

\author*[1]{\fnm{Louis} \sur{Callens}}\email{louis.callens@kuleuven.be}
\author[1]{\fnm{Bastiaan} \sur{Vandewal}}\email{bastiaan.vandewal@kuleuven.be}
\author[1]{\fnm{Ibrahim} \sur{Ibrahim}}\email{ibrahim.ibrahim@kuleuven.be}
\author[1]{\fnm{Jan} \sur{Swevers}}\email{jan.swevers@kuleuven.be}
\author[1]{\fnm{Wilm} \sur{Decré}}\email{wilm.decre@kuleuven.be}
\affil[1]{\orgdiv{MECO research team, department of Mechanical Engineering}, \orgname{KU Leuven}, \orgaddress{\street{Celestijnenlaan 300}, \city{Leuven}, \postcode{3001}, \country{Belgium}}}
\affil[1]{\orgdiv{Core Lab MPRO}, \orgname{Flanders Make@KU Leuven}, \country{Belgium}}

\abstract{
This paper proposes a novel and efficient optimization-based method for generating near time-optimal trajectories for holonomic vehicles navigating through complex but structured environments. The approach aims to solve the problem of motion planning for planar motion systems using magnetic levitation that can be used in assembly lines, automated laboratories or clean-rooms. In these applications, time-optimal trajectories that can be computed in real-time are required to increase productivity and allow the vehicles to be reactive if needed. The presented approach encodes the environment representation using free-space corridors and represents the motion of the vehicle through such a corridor using a motion primitive. These primitives are selected heuristically and define the trajectory with a limited number of degrees of freedom, which are determined in an optimization problem. As a result, the method achieves significantly lower computation times compared to the state-of-the-art, most notably solving a full Optimal Control Problem (OCP), OMG-tools or VP-STO without significantly compromising optimality within a fixed corridor sequence. The approach is benchmarked extensively in simulation and is validated on a real-world Beckhoff XPlanar system.
}

\keywords{Constrained motion planning, Optimization and Optimal Control, Motion Control, Collision Avoidance}

\maketitle

\section{Introduction}
\begin{figure}
    \centering
    \includegraphics[width=\linewidth, trim={0, 12.0cm, 0, 0cm}, clip]{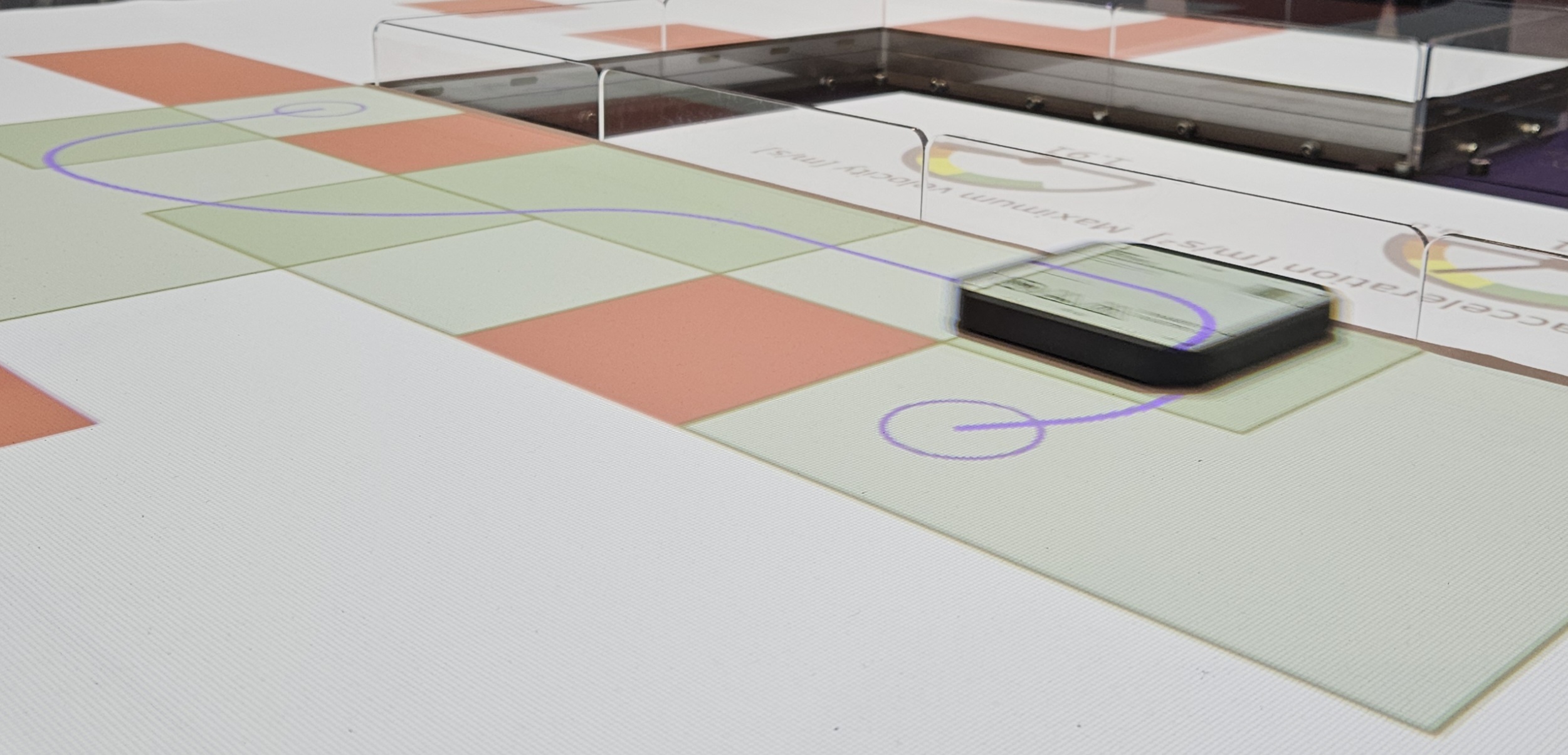}
    \caption{APM4220 XPlanar Mover (without any payload) by Beckhoff Automation, measuring 113mm $\times$ 113mm, executing a time-optimal trajectory that is projected from above as a blue line. The red squares represent (virtual) obstacles}
    \label{fig:hardware-xplanar}
\end{figure}

Mobile robots are utilized in industry in challenging environments where they are expected to autonomously navigate around obstacles. Holonomic vehicles, capable of moving in any direction independently and independent of their orientation, are especially appropriate for cluttered environments since they require less maneuvering space. Wheeled mobile robots can be made holonomic by using special wheel designs, turning them into Omnidirectional Mobile Robots (OMRs) \citep{omnidirectional-mobile-robots}. 
However, these omnidirectional wheels face challenges with slippage, uneven surfaces or wheel singularities. An alternative for wheeled holonomic vehicles are planar mover systems that use magnetic levitation to hover \citep{maglev-futuristic-overview}. In these systems, a table generating magnetic fields can move a permanent magnet around with micrometer precision \citep{planar-mover-system} \citep{planar-mover-system-2}. The high positioning accuracy and the contactless motion make this system exceptionally suited for small-scale assembly lines, cleanrooms or laboratories.
In recent years, this technology has become available in commercial products e.g. ACOPOS 6D by B\&R Industrial Automation, XBot Flyway by Planar Motor or the XPlanar Mover System by Beckhoff Automation (see Figure \ref{fig:hardware-xplanar}) and has already been adopted in industry to significantly increase flexibility in manufacturing. 

In order to fully exploit the flexibility offered by the hardware, a motion planning algorithm is required to route the movers from one point to another. There are a number of requirements for this algorithm. Firstly, it must be efficient in computation time such that it can be used to replan while moving to account for a dynamic environment. When the vehicles move fast, the planner must be especially efficient. Secondly, the motion planning algorithm must be able to adhere to velocity limits for safety and to acceleration limits to make sure there are no undesired sliding motions of an unsecured payload on the mover. Thirdly, the movers must be able to navigate through narrow passages, since the available free space will typically be limited due to high costs of the tiles that generate the magnetic fields. This means no conservatism can be introduced when considering the vehicles footprint. Finally, to increase the throughput of these motion systems in industrial applications, time-optimal motions are desired.

\subsection{Related work}
To tackle this motion planning problem, sampling-based planners could be used, such as Rapidly-exploring Random Trees (RRTs) \citep{RRT}. This approach was later extended to prove asymptotic optimality, leading to RRT$^*$ \citep{RRT*}. Attempts have been made to improve convergence. \citet{RRT-smart} proposed RRT-smart, which adds intelligent sampling techniques to sample close to obstacles and optimizes intermediate paths by changing the interconnections in the graph. \citet{informed-RRT} introduced informed-RRT, which samples the environment in a more direct way once an initial path has been found. An overview of RRT-based methods is given by \citet{RRT-survey}. It should be noted that these methods typically produce non-smooth, distance-optimal trajectories instead of time-optimal trajectories (with bounded acceleration). 

To account for vehicle dynamics within the sampling-based approach, \citet{kinodynamic-RRT} propose kinodynamic RRT$^*$ to compute optimal trajectories for vehicles with linear dynamics (such as a double integrator). Samples are connected with an optimized trajectory in closed-form instead of a straight line. However, the objective for these connecting trajectories contains a nonzero quadratic penalty on the acceleration in addition to the travel time, leading to trajectories that are not purely time-optimal. Acceleration limits are only implicitly considered by rejecting samples if bounds are exceeded, which can lead to many rejections and thus a large number of samples required when aiming for time-optimality. \citet{RRT-SQP-solver} extend the idea of Kinodynamic RRT$^*$ to nonlinear systems and address the issue of dealing with acceleration bounds indirectly. This is achieved by using Sequential Quadratic Programming (SQP) to solve the Boundary-Value-Problem (BPV) to find the constrained, time-optimal connections between sample points. However, optimality of the entire trajectory is still dependent on the samples themselves and is only achieved asymptotically. This leads to planning times that are too slow for use in real-time. Instead of optimizing the trajectories in between samples, Via-point-based Stochastic Trajectory Optimization (VP-STO) \citep{vp-sto} connects samples using the analytical solution of an OCP and optimizes the locations of the samples in a gradient-free way using Covariance Matrix Adaptation (CMA-ES), minimizing any user-specified objective. An online version of this algorithm is proposed as well. However, optimality is sacrificed by not running the optimization until convergence and by only using a limited number of via-points.

In contrast to sampling-based methods, optimization-based methods take into account dynamics, constraints and objective functions in a more direct way, by solving a so-called Optimal Control Problem (OCP).  However, solving a complete OCP is often too slow to be used in real-time, especially if many obstacles are present in the environment. To reduce the computation time for solving an OCP, \citet{omg1} proposed a spline-representation of the state trajectories and the controls in a toolbox called OMG-tools. This ensures smooth trajectories and reduces the degrees of freedom in the optimization problem, which allows it to be solved faster. However, the controls are not expected to be smooth in time-optimal trajectories. Obstacles are dealt with using separating hyperplanes which does not scale well if many obstacles are present. 

Instead of considering all obstacles in the optimization problem, the free-space can be modeled as a convex region to move through. One such approach is the Corridor Map Method \citep{corridor-map-method,corridor-map-method-improved}. Large, elliptical, obstacle-free regions can be constructed quickly using Iterative Regional Inflation by Semidefinite programming (IRIS) \citep{large-convex-regions}. This algorithm alternates between the computation of separating hyperplanes and the inflation of the ellipse inscribed in the polytope defined by these hyperplanes. After constructing an initial guess using a search algorithm, \citet{iwanaga} use a modified version of the IRIS algorithm to produce polytope corridors around the OCP grid points. This approach allows to construct optimal trajectories through cluttered environments with tight passages for a bycicle model. OMG-tools was extended to use basic rectangular corridors representing the free-space to plan an optimal trajectory for a holonomic vehicle \citep{omg2}. Rectangular corridors can represent the free-space accurately in structured environments. \citet{truck-trailer} applied the idea of using rectangular corridors to parking maneuvers for a truck with a trailer. Here, corridors represent lanes on the road, parking lots or allowed maneuvering spaces. \citet{sonia} exploit the simplification of the environment represented by using only rectangular corridors to determine the time-optimal trajectory analytically, eliminating the need for optimization completely. This is done by selecting exact motion primitives based on the geometry of the rectangular corridors. However, this planner is limited to velocity-control of a unicycle and therefore produces trajectories with infinite acceleration, which is undesired in the case of a planar mover system carrying payloads. An extension to selecting optimal primitives that include acceleration limits analytically  is unclear.

Both OMG-tools and \citet{truck-trailer}, use IPOPT \citep{ipopt} to solve the resulting optimization problems, an interior-point-based general purpose solver for Nonlinear Programs (NLPs). However, IPOPT does not exploit the sparsity structure in the Lagrangian Hessian and constraint Jacobian that typically arises in OCPs. FATROP \citep{vanroye2023fatrop} is based on the ideas of IPOPT, but does exploit this known structure, significantly reducing the computation time to solve the OCP. Both IPOPT and FATROP can easily be used from \texttt{CasADi} \citep{casadi}, an open-source toolbox that constructs symbolic expressions that can easily be derived using Algorithmic Differentiation (AD). The Opti stack, part of \texttt{CasADi}, allows to easily specify optimization problems and solve them with a solver of choice.

\subsection{Contributions}
In this paper, we present an optimization-based approach to efficiently compute collision-free and near time-optimal trajectories through structured environments for the holonomic vehicle, while adhering to velocity and acceleration bounds. Our approach models the free space in the structured environment using rectangular corridors. Similar to the ideas of \cite{vp-sto} and \cite{omg2}, a low dimensional representation of trajectory is obtained and optimized. We achieve this by selecting a motion primitives based on the geometry of the corridors, related to the work by \citet{sonia}.

We extensively benchmarked the approach and validated both the selection and optimization of the motion primitives given a corridor sequence, showing that in most cases, we find the optimal solution in significantly less time. In addition, the idea of planning through corridor sequences is evaluated by benchmarking against VP-STO. Apart from benchmarking in simulation, the approach is validated on a real system: the XPlanar Mover System by Bechhoff Automation. This validation demonstrates that the proposed trajectories can accurately be executed on real hardware.

The main contributions of this work are (1) the selection of time-optimal motion primitives based only on the geometry of rectangular corridors, and (2) the formulation of the resulting optimization problem, which can be solved efficiently without sacrificing optimality significantly to determine the remaining degrees of freedom.

\subsection{Structure of the paper}
The remainder of this paper is structured as follows, we will first detail the problem description in Section \ref{sec:problem-description} before explaining the proposed method in Section \ref{sec:methodology}. The proposed method is benchmarked and analyzed in Section 
\ref{sec:benchmark}. Validation on a real-system is discussed in Section \ref{Sec:experimental-validation}. Finally, we formulate conclusions and suggestions for future work in Section \ref{sec:conclusions}.

\subsection{Note on Notation}
In this paper, bold lower case letters represent two-dimensional vectors unless specified otherwise. For example $\bm{q} = (q^x, q^y)^\top \in \mathbb{R}^2$. Operations on these vectors are always applied element-wise.
\section{Problem Description} \label{sec:problem-description}
The dynamics considered in this work are those of a holonomic vehicle, modeled as a two-dimensional double integrator. We consider a vehicle with a rectangular footprint of width $W$ and length $L$. The vehicle is aligned with the world frame axes and does not rotate. Note that the XPlanar mover is only able to translate when the vehicle is aligned with the axes, justifying the assumptions made on the vehicle's orientation. The dynamics are given by
\begin{align}
    \dot{\bm{p}}(t) &= \bm{v}(t) &
    \dot{\bm{v}}(t) &= \bm{a}(t)
\end{align}
where $\bm{p}(t)$, $\bm{v}(t)$, $\bm{a}(t)$ represent the geometrical center, the velocity and the acceleration of the vehicle respectively.

The goal is to move the vehicle from a starting position $\bm{p}_0$ to a destination $\bm{p}_n$ as fast as possible with some given initial velocity $\bm{\bar{v}}_0$ and zero terminal velocity $\bm{v}_n=0$, while adhering to velocity and acceleration limits given by
\begin{align}
    ||\bm{v}(t)||_\infty &\leq v_{\textrm{max}} & 
    ||\bm{a}(t)||_\infty &\leq a_{\textrm{max}}.\label{eq:velocity-acceleration-limits}
\end{align}
The use of the infinity-norm models holonomic vehicles with independent actuation in each dimension. Holonomic vehicles such as gantry cranes clearly have decoupled actuation. For planar mover systems, maximal acceleration values are more complicated and depend on center of mass of the payload, height of the mover and other factors. The proposed methodology uses the infinity-norm to decouple and hence simplify the planning problem which reduces computation time. 

The vehicle must avoid obstacles present in the environment. We assume that the environment is static and that an occupancy grid map of this environment is available. The grid is a matrix aligned with the world axes where each grid cell has dimensions at least $W \times L$. The proposed method targets structured environments that can be represented in this way. For planar motion systems using magnetic levitation, the environment typically consists of connected square tiles that generate the magnetic fields, meaning the realworld environment naturally fits the proposed environment representation. Note that in this application, obstacles will typically be stationary movers, environment boundaries or faulty tiles. Parking lots or warehouses are other examples of environments that can easily be represented this way.

Since the dynamics at hand are decoupled, point-to-point motion can be computed separately for both dimensions. As such, there will always be one dimension that acts as a bottleneck when moving in minimum-time. Suppose the optimal moving times are computed for both dimensions and $T^*_x < T^*_y$, then any feasible solution in the $x$ direction that satisfies $T_x < T^*_y$ has the same overall objective function value. This means there is additional freedom in the $x$ dimension. In this case, the $y$ dimension can be considered the bottleneck direction ($B$) whereas the $x$ dimension is the free direction ($F$).
\section{Methodology} \label{sec:methodology}
\begin{figure*}[!h]
    \centering
    \begin{subfigure}[t]{0.24\linewidth}
        \includegraphics[width=\linewidth]{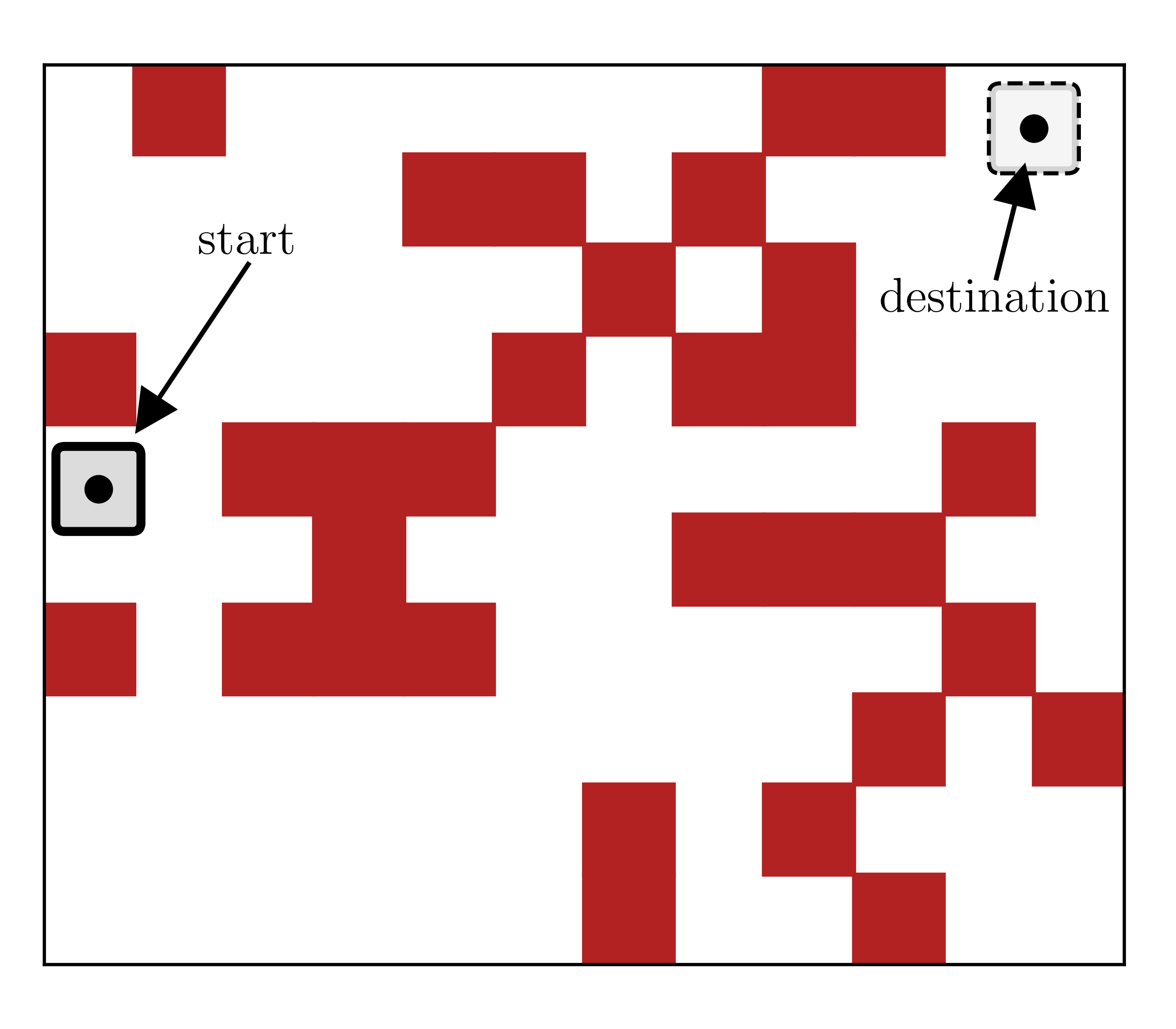}
        \subcaption{}
        \label{fig:corridor-construction-env}
    \end{subfigure}
    \begin{subfigure}[t]{0.24\linewidth}
        \includegraphics[width=\linewidth]{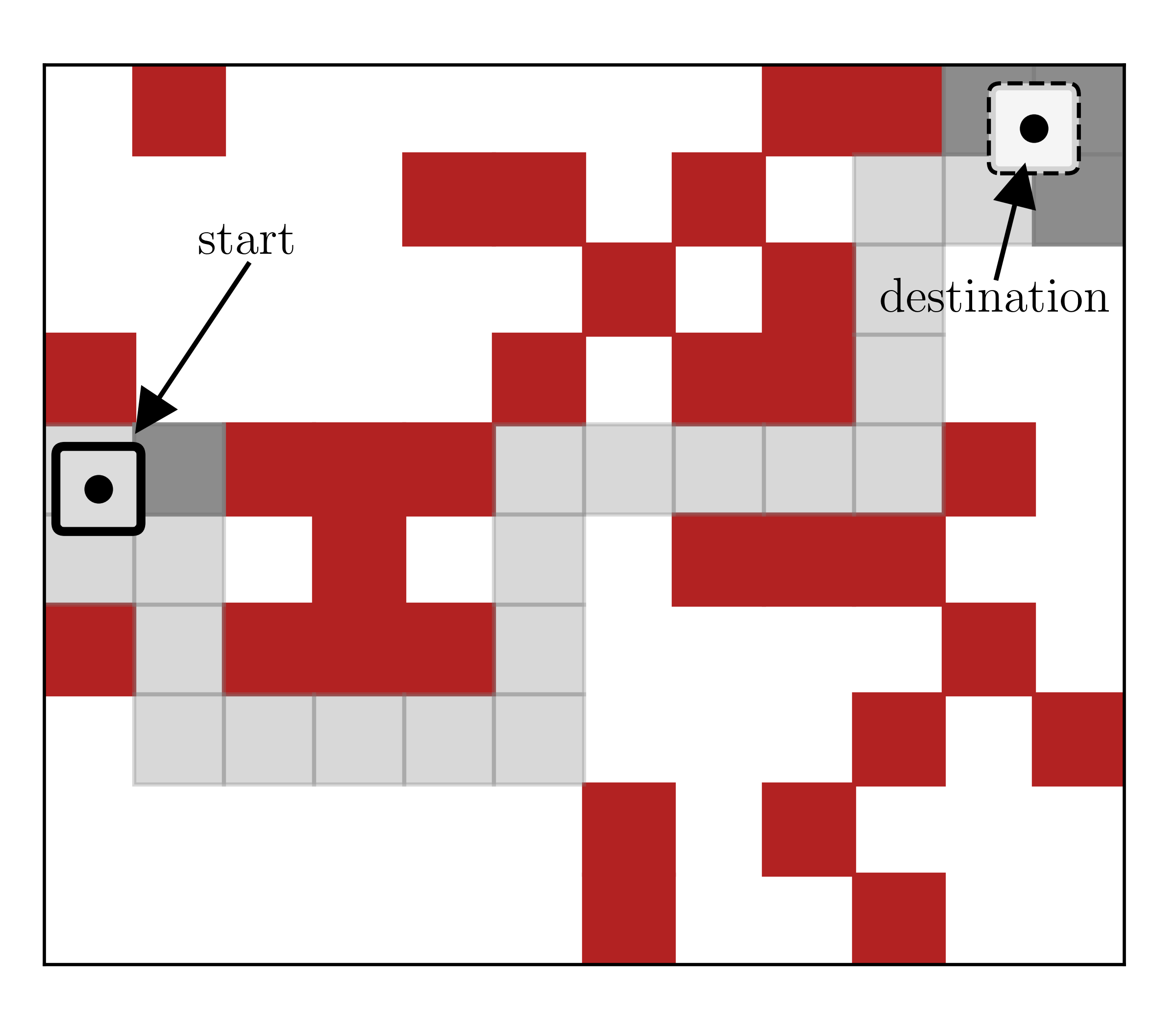}
        \subcaption{}
        \label{fig:corridor-construction-path}
    \end{subfigure}
    \begin{subfigure}[t]{0.24\linewidth}
        \includegraphics[width=\linewidth]{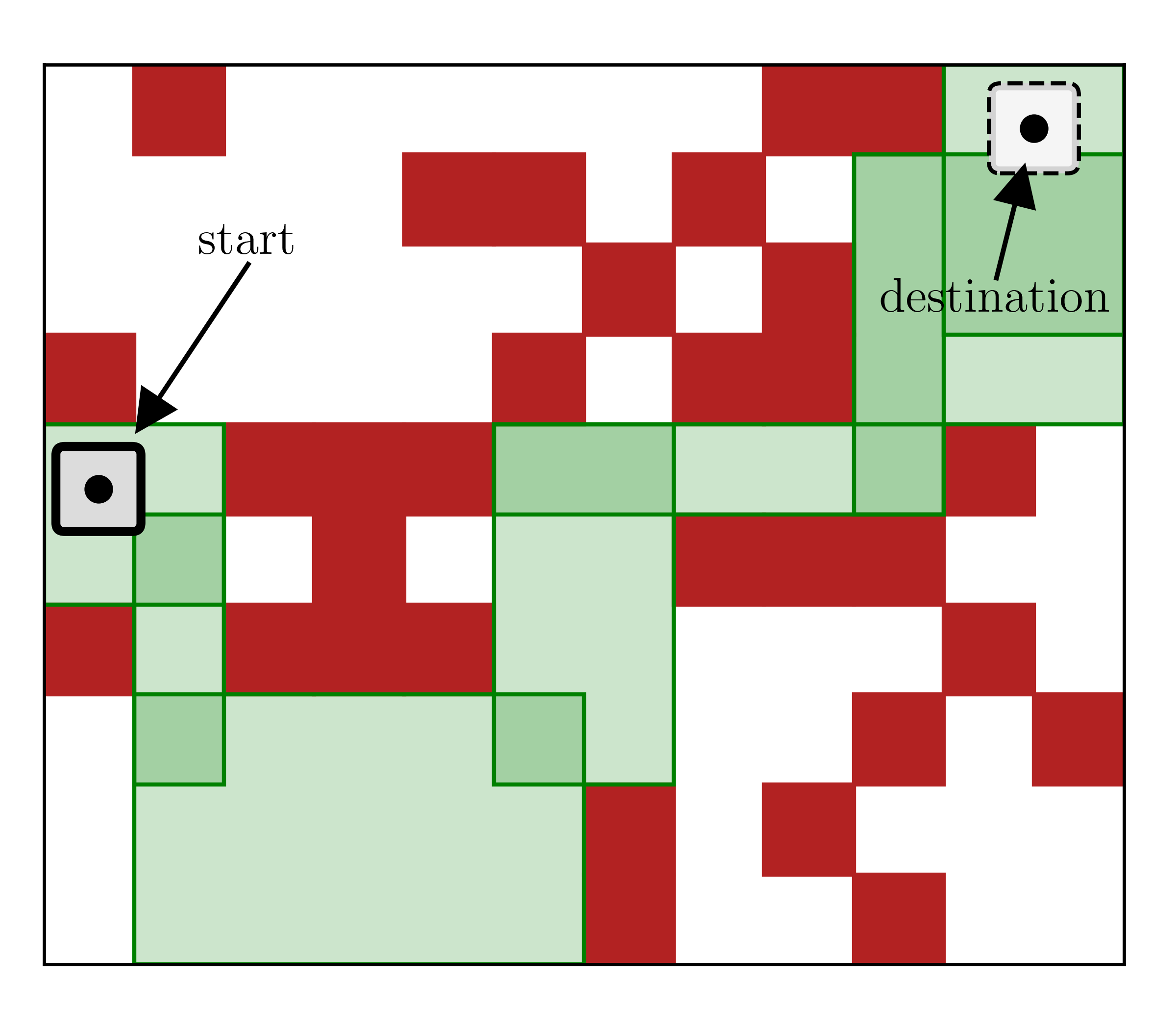}
        \subcaption{}
        \label{fig:corridor-construction-final}
    \end{subfigure}
    \begin{subfigure}[t]{0.24\linewidth}
        \includegraphics[width=\linewidth]{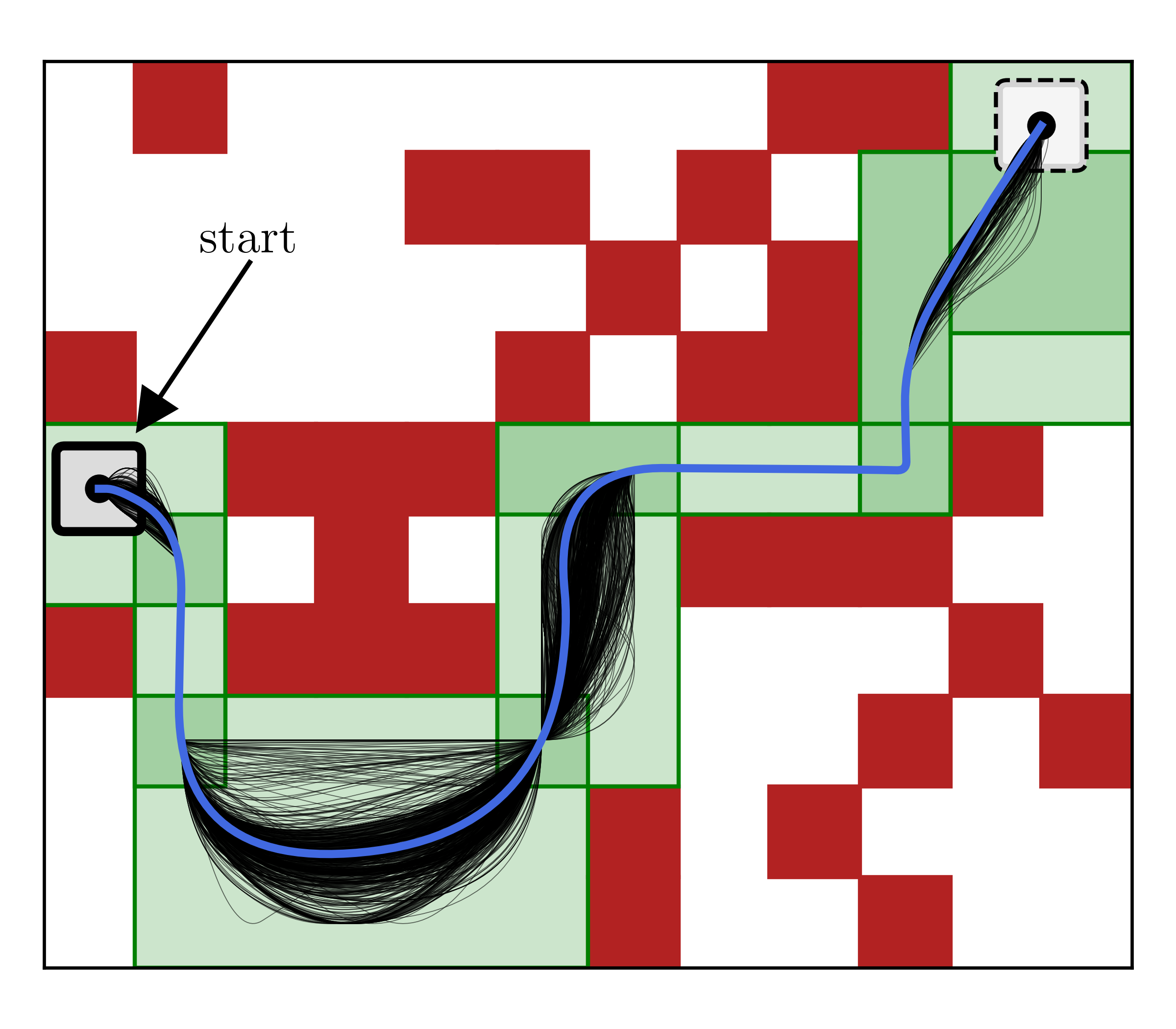}
        \subcaption{}
        \label{fig:corridor-construction-trajs}
    \end{subfigure}
    \caption{Illustration of the different steps of the proposed method. The red squares indicate cells that are occupied by an obstacle. The starting and final positions are shown as well as the vehicle footprint. (a) The environment with starting position and desired destination indicated. (b) The path found by performing a graph-search, denoted by $\mathcal{P}$, is shown in light gray color. The cells in $\mathcal{P}'$ that are not in $\mathcal{P}$ are shown in dark gray. (c) The corridor sequence $\bm{C}$ is shown in light green. The overlap of two consecutive corridors is shown in dark green. (d) Black lines show 500 trajectories described using the selected motion primitives (random samples from $\bm{\Omega}_n$). The blue trajectory is the time-optimal one ($\Omega_n^* = \bm{\Phi}\left(\Pi^*\right)$ where $\Pi^* \in \mathcal{X}$)}
    \label{fig:corridor-construction}
\end{figure*}

The proposed method for solving the problem defined above consists of three steps and is illustrated in Figure \ref{fig:corridor-construction}. First, a sequence of rectangular corridors is constructed to simplify the representation of the environment. This is explained in Section \ref{sec:corridors} and visualized in Figures \ref{fig:corridor-construction-path}-\ref{fig:corridor-construction-final}. Second, a dynamically feasible and optimal trajectory is planned through these corridors. As shown in Section \ref{sec:analytical}, an analytical solution ignoring corridor bounds is constructed. If this analytical solution is collision-free, the problem is solved. If not, the corridor bounds have to be taken into account explicitly. This is achieved by selecting parametric motion primitives (defined in Section \ref{sec:primitive}) based on the geometry of the corridors (Section \ref{sec:primitive-selection}). Figure \ref{fig:corridor-construction-trajs} shows random trajectories represented by the selected primitives for an example. The remaining degrees of freedom are optimized as explained in Section \ref{sec:optimization}. Some specifics of the formulation of the optimization problem are discussed in Section \ref{sec:details}. Section \ref{sec:overview} provides a summary of the approach in pseudocode. Finally, limitations of the method are listed in Section \ref{sec:limitations}.

Since the proposed method is based on parametric motion primitives, we refer to it as PMP. An implementation is available at \url{https://github.com/meco-group/parametric-motion-planner#}.

\subsection{Construction of a Sequence of Free-Space Corridors} \label{sec:corridors}
The goal of constructing a sequence of $n$ free-space corridors $\bm{C}$ is to represent the environment in a simplified and convex way. The corridors are rectangles that are aligned with the world axes, so they can be defined as
\begin{equation}
    \mathcal{C}_i \coloneq \{\bm{q} : q^x \in [x^i_{\mathrm{min}}, x^i_{\mathrm{max}}], q^y \in [y^i_{\mathrm{min}}, y^i_{\mathrm{max}}]\}
\end{equation}
where $i = 0,\ldots,n-1$. This definition and choice of rectangular corridors allows to treat the constraints in a decoupled way, which is not the case if elliptical corridors are considered such as those computed by IRIS \citep{large-convex-regions}. More precisely, enforcing the vehicle at position $\bm{p}$ to be within $\mathcal{C}_i$ leads to the constraints
\begin{align}
    x_{\mathrm{min}}^i + \frac{W}{2} \leq & p_x \leq x_{\mathrm{max}}^i - \frac{W}{2}\label{eq:corridor-bound-x}\\
    y_{\mathrm{min}}^i + \frac{L}{2} \leq & p_y \leq y_{\mathrm{max}}^i - \frac{L}{2}\label{eq:corridor-bound-y}.
\end{align}

Define $\mathcal{O}_{i,j} \coloneq \mathcal{C}_{i} \cap \mathcal{C}_j$ as the intersection of corridors $i$ and $j$. We demand $\mathcal{O}_{i, j} = \emptyset$ if $|i - j| > 1$ and $\mathcal{O}_{i-1, i}$ large enough to contain the vehicle $\forall i \in \{1,\ldots,n-1\}$.

The first step in the construction of $\bm{C}$ is to compute a path as a sequence of grid cells that are denoted $g_i$ that is defined by
\begin{equation}
\begin{aligned}
    \mathcal{P} \coloneq (&g_0 \dotsi g_{m} : g_i \text{ is free, }\\
    &D_M(g_i, g_{i+1}) = 1, \\
    &\bm{p}_0 \in g_0, \bm{p}_n \in g_{m})
\end{aligned}
\end{equation}
where $D_M$ represents the Manhattan distance in the occupancy grid. $\mathcal{P}$ is computed by performing a breadth-first search. Define the sets
\begin{align}
    \mathcal{F}(\bm{p}) &\coloneq \left\{\bm{q} : \left|\bm{p} - \bm{q}\right| \leq \begin{bmatrix} \frac{W}{2} & \frac{L}{2} \end{bmatrix}^\top \right\}\label{eq:set-footprint},\\
    \mathcal{S}(\bm{p}) &\coloneq \{g_i : \exists \bm{q} \in \mathcal{F}(\bm{p}) \cap g_i\}.
\end{align}
$\mathcal{F}$ represents the footprint of the vehicle and contains all positions occupied by the vehicle at position $\bm{p}$. $\mathcal{S}$ is the set of cells that are (partially) occupied by the vehicle at position $\bm{p}$. To ensure the vehicle can satisfy equations \eqref{eq:corridor-bound-x}-\eqref{eq:corridor-bound-y} at $\bm{p}_0$ and $\bm{p}_n$, $\mathcal{F}(\bm{p}_0) \subseteq \mathcal{C}_0$ and $\mathcal{F}(\bm{p}_n) \subseteq \mathcal{C}_{n-1}$ must hold. To this end, we extend $\mathcal{P}$ with $\mathcal{S}(\bm{p}_0)$ and $\mathcal{S}(\bm{p}_n)$. This leads to the definition
\begin{equation}
    \mathcal{P}' \coloneq \mathcal{S}(\bm{p}_0) \smallfrown \mathcal{P} \smallfrown \mathcal{S}(\bm{p}_n)
\end{equation}
where $\smallfrown$ denotes concatenation. $\mathcal{P}$ and $\mathcal{P}'$ are shown for an example environment in Figure \ref{fig:corridor-construction-path}. The darker cells are added such that the complete vehicle footprint at $\bm{p}_0$ and $\bm{p}_n$ is included. Because of the assumption on the dimension of the grid cells, some feasible path from $\bm{p}_0$ to 
$\bm{p}_n$ exists in $\mathcal{P}'$.

\begin{figure}
    \centering
    \begin{subfigure}{0.45\linewidth}
        \includegraphics[width=\linewidth]{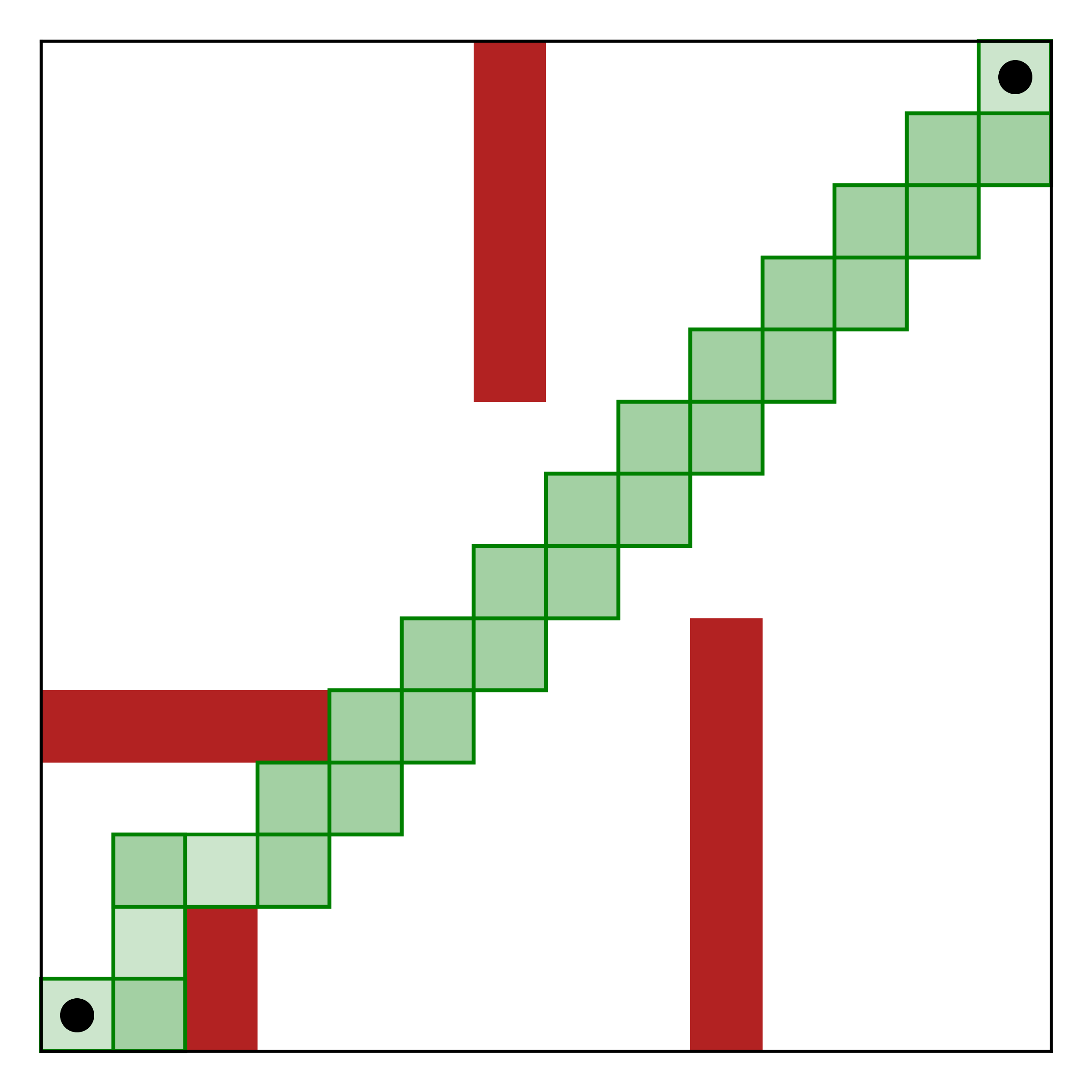}
        \subcaption{initial $\bm{C}$}
    \end{subfigure}
    \begin{subfigure}{0.45\linewidth}
        \includegraphics[width=\linewidth]{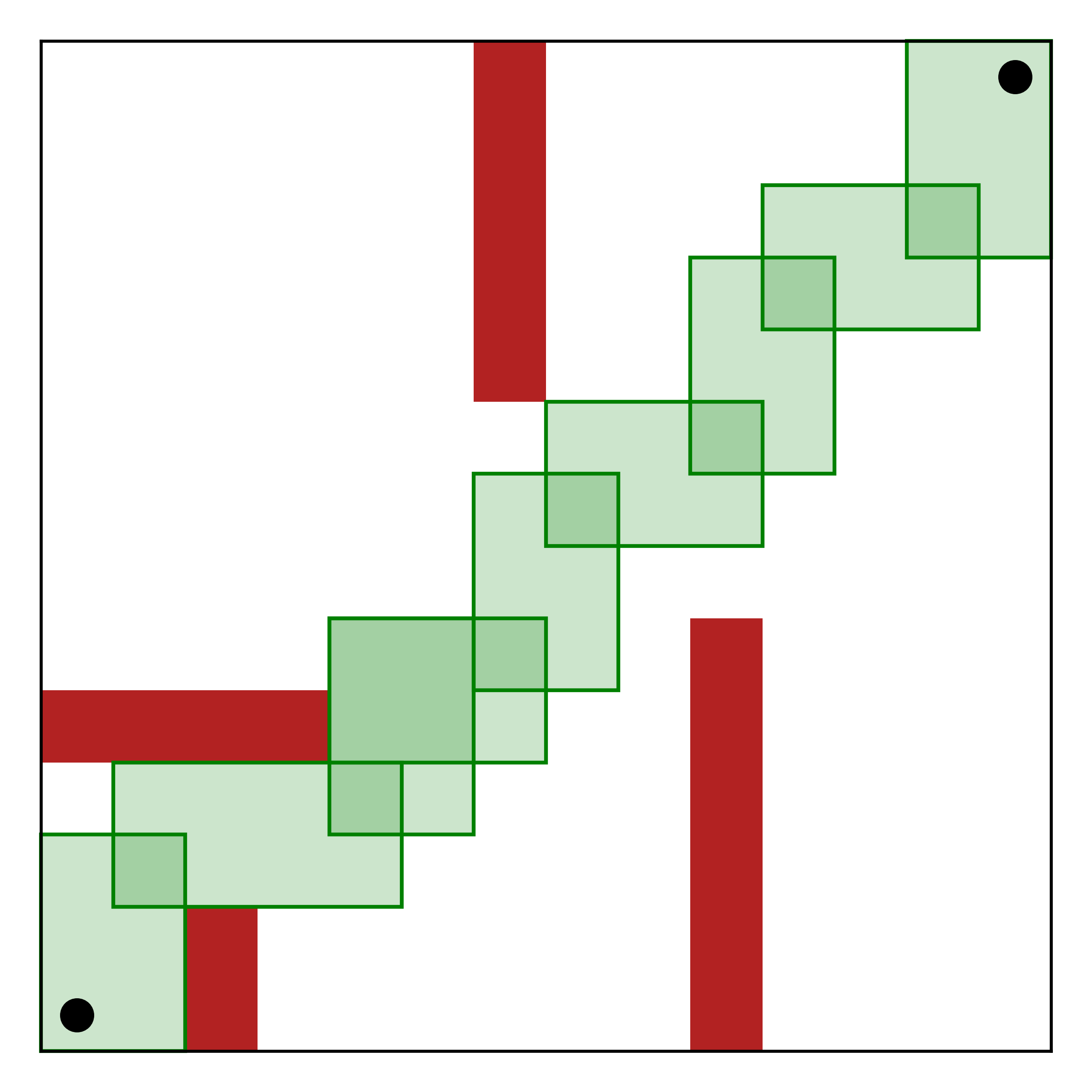}
        \subcaption{$\bm{C}$ after 1 iteration}
    \end{subfigure}
    \begin{subfigure}{0.45\linewidth}
        \includegraphics[width=\linewidth]{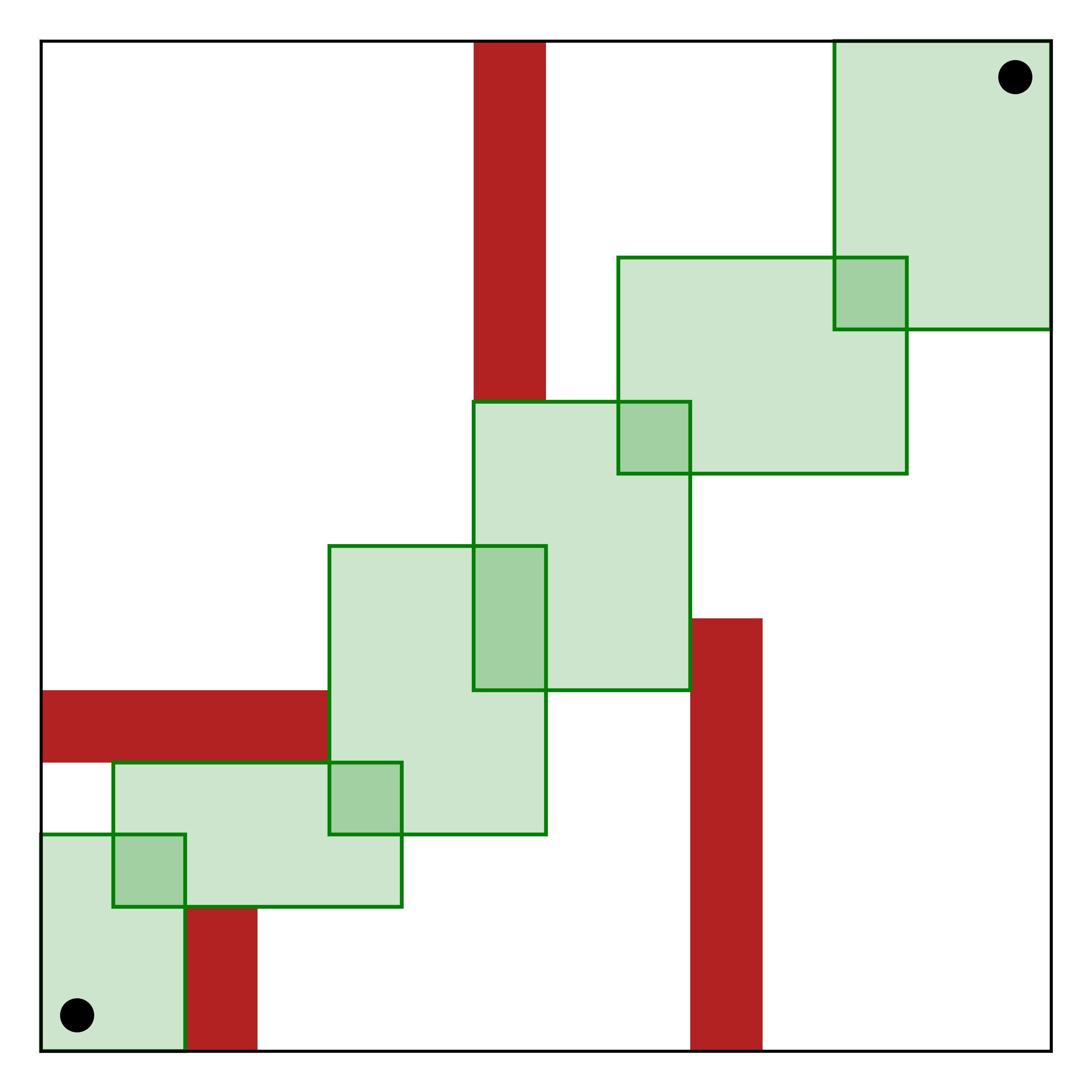}
        \subcaption{$\bm{C}$ after 2 iterations}
    \end{subfigure}
    \begin{subfigure}{0.45\linewidth}
        \includegraphics[width=\linewidth]{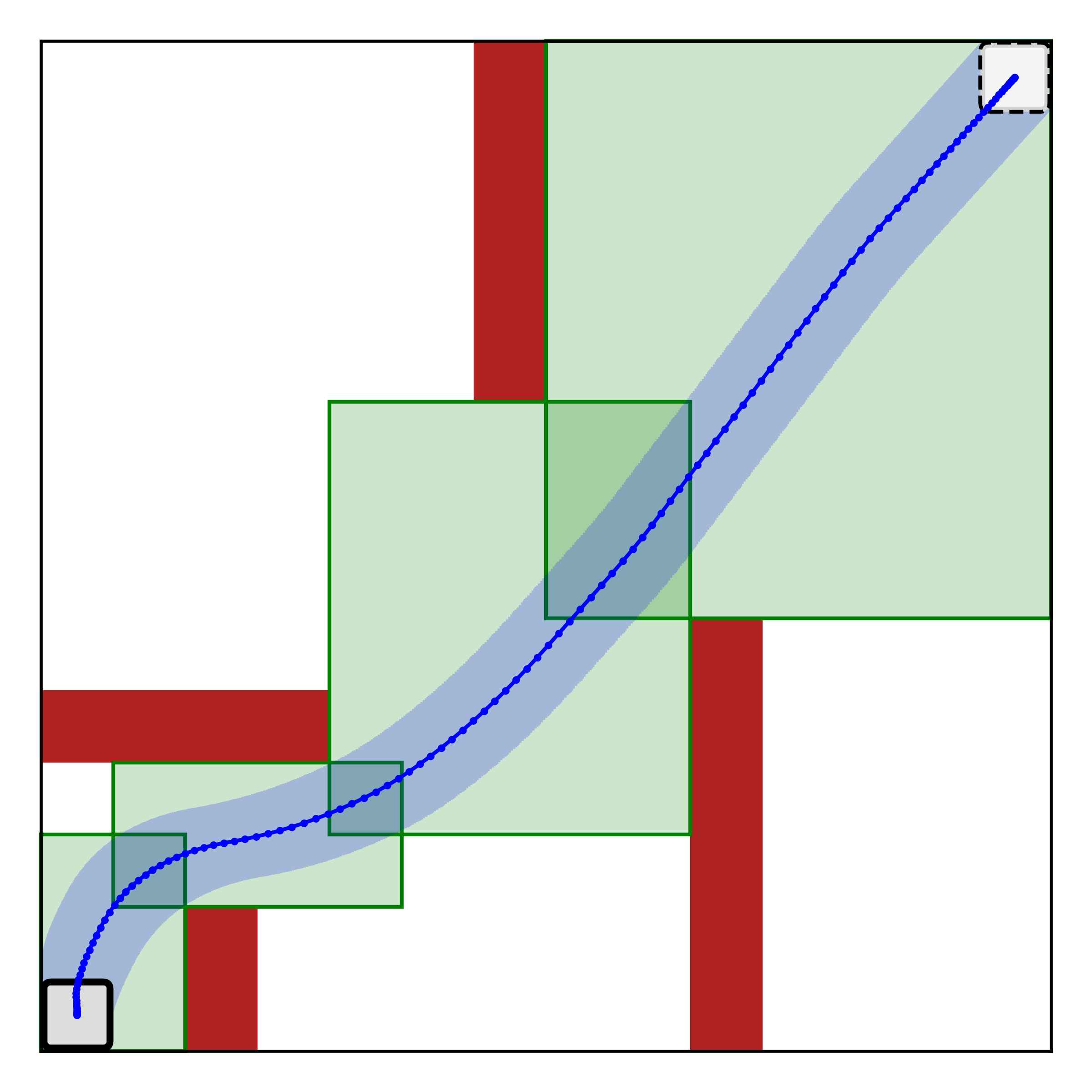}
        \subcaption{final $\bm{C}$ (7 iterations)}
    \end{subfigure}
    \caption{Illustration of the iterative growing strategy for an example environment. Note that corridors $\mathcal{C}_i$ are removed if $\mathcal{C}_i \subseteq \left(\mathcal{C}_{i-1} \cup \mathcal{C}_{i+1} \right)$ or if $\mathcal{C}_{i-1} \cap \mathcal{C}_{i+1}$ is nonempty, which occurs often in this environment. The corridor sequence $\bm{C}$ is shown in light green. The overlap of two consecutive corridors is shown in dark green. Obstacles are shown in red. The blue line shows the optimized trajectory.}
    \label{fig:corridor-growing}
\end{figure}

$\bm{C}$ is obtained by splitting $\mathcal{P}'$ in a series of contiguous intervals of cells that all lie on the same row or column.
To ensure each corridor to be as wide as possible allowing the vehicle as much free space as possible, an iterative growing strategy is applied to $\bm{C}$ by checking occupancy of the grid cells along the sides of each corridor. After every enlargement iteration, corridors $\mathcal{C}_i$ are removed from the sequence if $\mathcal{C}_i \subseteq \left(\mathcal{C}_{i-1} \cup \mathcal{C}_{i+1} \right)$ or if $\mathcal{C}_{i-1} \cap \mathcal{C}_{i+1}$ is nonempty. This process is illustrated in Figure \ref{fig:corridor-growing}. Figure \ref{fig:corridor-construction-final} shows the constructed corridors for an example. 

\subsection{Analytical Solution without Considering Obstacles} \label{sec:analytical}
Suppose that the corridors are infinitely large so any trajectory is collision-free. Moreover, suppose that problem is solved in each dimension separately, resulting in optimal moving times $T^*_x$ and $T^*_y$ and suppose, without loss of generality, that $T^*_x < T^*_y$. Choose $a_x(t) = 0 \quad \forall t \in (T^*_x, T^*_y]$. Because $v_x(T^*_x) = 0$, this ensures that $p_x(T^*_y) = p_n^x$ and $v_x(T^*_y) = 0$. In other words, the movement in both dimensions can be completely decoupled because the terminal velocity is considered to be zero.

The decoupling of the two directions reduces the problem to the constrained minimum-time motion of a double integrator, which can be solved completely using Pontryagin's Maximum Principle, as described by \citet{constrained-double-integrator}.
If all points on the resulting trajectory satisfy \eqref{eq:corridor-bound-x} and \eqref{eq:corridor-bound-y} for some corridor in $\bm{C}$, a feasible solution has been found. If not, parametric primitives are selected and optimized as described in the next sections.

\subsection{Definition of the Parametric Motion Primitive} \label{sec:primitive}
The time-optimal acceleration profile for a point-to-point motion in one dimension is generally given by 
\begin{equation}
    a(t) = 
    \begin{cases} 
        \alpha_k a_{\textrm{max}} & t \in [t_0, t_1)\\
        0 & t \in [t_1, t_2)\\
        \alpha_{k+1} a_{\textrm{max}} & t \in [t_2, t_3)        
    \end{cases}
\end{equation}
where $\alpha_i \in \{-1, 1\}$ and can be referred to as a bang-coast-bang acceleration.

Let $\tau_i \coloneq t_{i+1} - t_i$ be the time-durations of each acceleration phase. Integrating this acceleration profile leads to a velocity profile consisting of a linear segment, a constant segment and another linear segment.

\begin{figure}
    \centering
    \includegraphics[width=\linewidth]{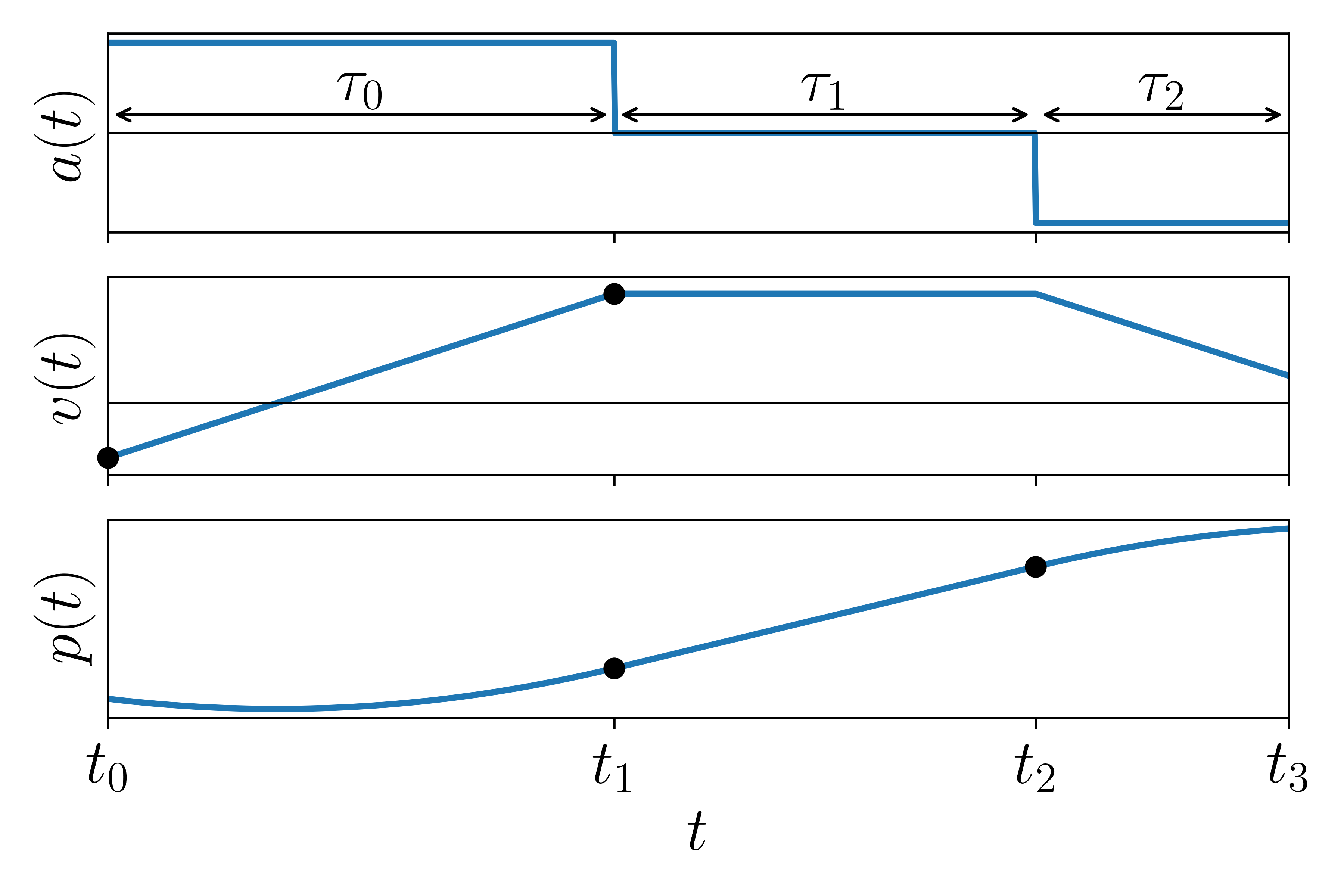}
    \caption{Illustration of the one-dimensional motion primitive $\phi(1, -1, p_k, v_k, \bm{\tau})$. The black dots are points on which box constraints are applied, as mentioned in Section \ref{sec:optimization}}
    \label{fig:primtive-illustration}
\end{figure}

This, in turn, leads to a position profile which has a parabolic segment, a linear segment followed by another parabolic segment. This one-dimensional movement can be described fully by the primitive
\begin{equation}
    \phi\left(\alpha_k, \alpha_{k+1}, p_k, v_k, \bm{\tau} \right) \coloneq \begin{bmatrix} p(t) & v(t) & a(t) \end{bmatrix}^\top
\end{equation} where $p(t_0) = p_k$, $v(t_0) = v_k$ and $\bm{\tau} \coloneq \begin{bmatrix} \tau_0 & \tau_1 & \tau_2 \end{bmatrix}$ is the vector collecting the time-durations of the accelerations for a single primitive in one dimension. 
An example is shown in Figure~\ref{fig:primtive-illustration}.

The two-dimensional motion primitive is then given by
{\small
\begin{equation*}
    \Phi\left(\bm{\alpha}_k, \bm{\alpha}_{k+1}, \bm{p}_k, \bm{v}_k, \bm{\tau}^x, \bm{\tau}^y \right) \coloneq \begin{bmatrix}
        \phi\left(\alpha_k^x, \alpha_{k+1}^x, p_k^x, v_k^x, \bm{\tau}^x \right) \\[1.5ex]
        \phi\left(\alpha_k^y, \alpha_{k+1}^y, p_k^y, v_k^y, \bm{\tau}^y \right)
    \end{bmatrix}
\end{equation*}
}
Let $T^{(.)} \coloneq \sum_{i=0}^{2}{\tau_i^{(.)}}$ be the total time of the primitive in one dimension, then $T \coloneq T^x = T^y$ must hold. This constrained is required to ensure that the primitive reaches its end-point at the same time in each dimension. Suppose we choose the parameters such that $p^x(T^x) = \tilde{p}^x$ and $p^y(T^y) = \tilde{p}^y$. In that case, the trajectory of the vehicle will only actually pass through the point $(\tilde{p}^x, \tilde{p}^y)$ if the timing constraint holds. In addition, the primitive would not be well defined for $t \in (\mathrm{min}(T^x, T^y), \mathrm{max}(T^x, T^y)]$. Let $\bm{\lambda}_i \coloneq \begin{bmatrix} \tau_i^x & \tau_i^y \end{bmatrix}^\top$. Now we can define the function
\begin{align*}
    f\left(\Phi\right) &\coloneq
    \begin{bmatrix}
        \hspace{-2em}\bm{p}_k + \bm{v}_k T + a_{\textrm{max}} \bm{\alpha}_k \bm{\lambda}_0 \left(\bm{\lambda}_1 + \bm{\lambda}_2\right) \\
        \hspace*{4em}+ 0.5 a_{\textrm{max}} \left( \bm{\alpha}_k \bm{\lambda}_0^2 + \bm{\alpha}_{k+1} \bm{\lambda}_2^2 \right) \\[1.5ex]
        \bm{v}_k + a_{\textrm{max}}\bm{\alpha}_k \bm{\lambda}_1 + a_{\textrm{max}}\bm{\alpha}_{k+1} \bm{\lambda}_3
    \end{bmatrix}
\end{align*}
which analytically integrates the position and velocity of the primitive.

The complete trajectory is described by a sequence of $n$ parametric primitives: one primitive $\Phi$ for each corridor $\mathcal{C}_k$.
In Section \ref{sec:primitive-selection}, a strategy is explained to determine $\{\bm{p}_k\}_{k=1}^{n-1}$ and $\{\bm{\alpha}_k\}_{k=0}^n$ so we can assume these to be given. Only the durations $\bm{\tau}_k^x$ and $\bm{\tau}_k^y$ and starting velocities $\bm{v}_k$ for every primitive are left to be determined. This means there are $8n$ remaining parameters that can be collected in a vector $\Pi \in \mathbb{R}^{2n} \times \mathbb{R}_{\geq0}^{6n}$ as
\begin{equation}
    \Pi \coloneq \begin{bmatrix} \bm{v}_0 \dotsi \bm{v}_{n-1} & \bm{\tau}_0^x \dotsi \bm{\tau}_{n-1}^x & \bm{\tau}_0^y \dotsi \bm{\tau}_{n-1}^y \end{bmatrix}
\end{equation}
which completely determines the set of $n$ primitives 
\begin{equation}
    \bm{\Phi}_n(\Pi)~\coloneq~\{\Phi_k\left(\bm{\alpha}_k, \bm{\alpha}_{k+1}, \bm{p}_k, \bm{v}_k, \bm{\tau}_k^x, \bm{\tau}_k^y \right)\}_{k=0}^{n-1}.
\end{equation}

Some degrees of freedom are eliminated by applying necessary constraints on the trajectory. More precisely, every primitive has a timing constraint enforcing the sum of the time durations in both dimensions to be equal as described earlier ($n$ constraints). On top of that, every primitive must reach a known final position and either ensure continuity of the velocity or reach the terminal velocity ($4n$ constraints). Finally, an initial velocity is given ($2$ constraints). This leads to $3n-2$ effective degrees of freedom and the definition of the feasible set of parameters as
\begin{equation}
\begin{aligned}
    \mathcal{X} \coloneq \{&\Pi \in \mathbb{R}^{2n} \times \mathbb{R}_{\geq0}^{6n} : \hspace{0.5em}\bm{v}_0 = \bm{\bar{v}}_0,\hspace{0.5em}\\
    &\!\begin{aligned}[t]
        &\textstyle\sum_{i=0}^{2}{\tau_{k,i}^x} = \sum_{i=0}^{2}{\tau_{k,i}^y} & &\hspace{0.5em} k = 0...n-1,\\
        &f(\Phi) = \begin{bmatrix} \bm{p}_{k+1} & \bm{v}_{k+1}\end{bmatrix}^\top & &\hspace{0.5em} k = 0...n-1\}.
    \end{aligned}
\end{aligned}
\end{equation}
Note that we do not explicitly eliminate parameters in $\Pi$ using these constraints. Instead, all parameters are kept as optimization variables and are subject to the above constraints, as shown in equation \ref{eq:opti-feasibility}.

We consider the set of trajectories $\bm{\Omega}_n$ defined as
\begin{equation}
    \bm{\Omega}_n \coloneq \{\bm{\Phi}_n(\Pi) : \Pi \in \mathcal{X}\}.
\end{equation}

The black lines in Figure \ref{fig:corridor-construction-trajs} show trajectories sampled from $\bm{\Omega}_n$ for an example environment.

\subsection{Selection of the Parametric Motion Primitives} \label{sec:primitive-selection}
Based on the geometry of the corridor sequence $\bm{C}$, the values $\{\bm{\alpha}_k\}_{k=0}^n$ and the points $\{\bm{p}_k\}_{k=1}^{n-1}$ can be determined heuristically. $\bm{p}_0$ and $\bm{p}_n$ follow from the problem definition. 

\begin{figure}
    \centering
    \includegraphics[width=0.8\linewidth]{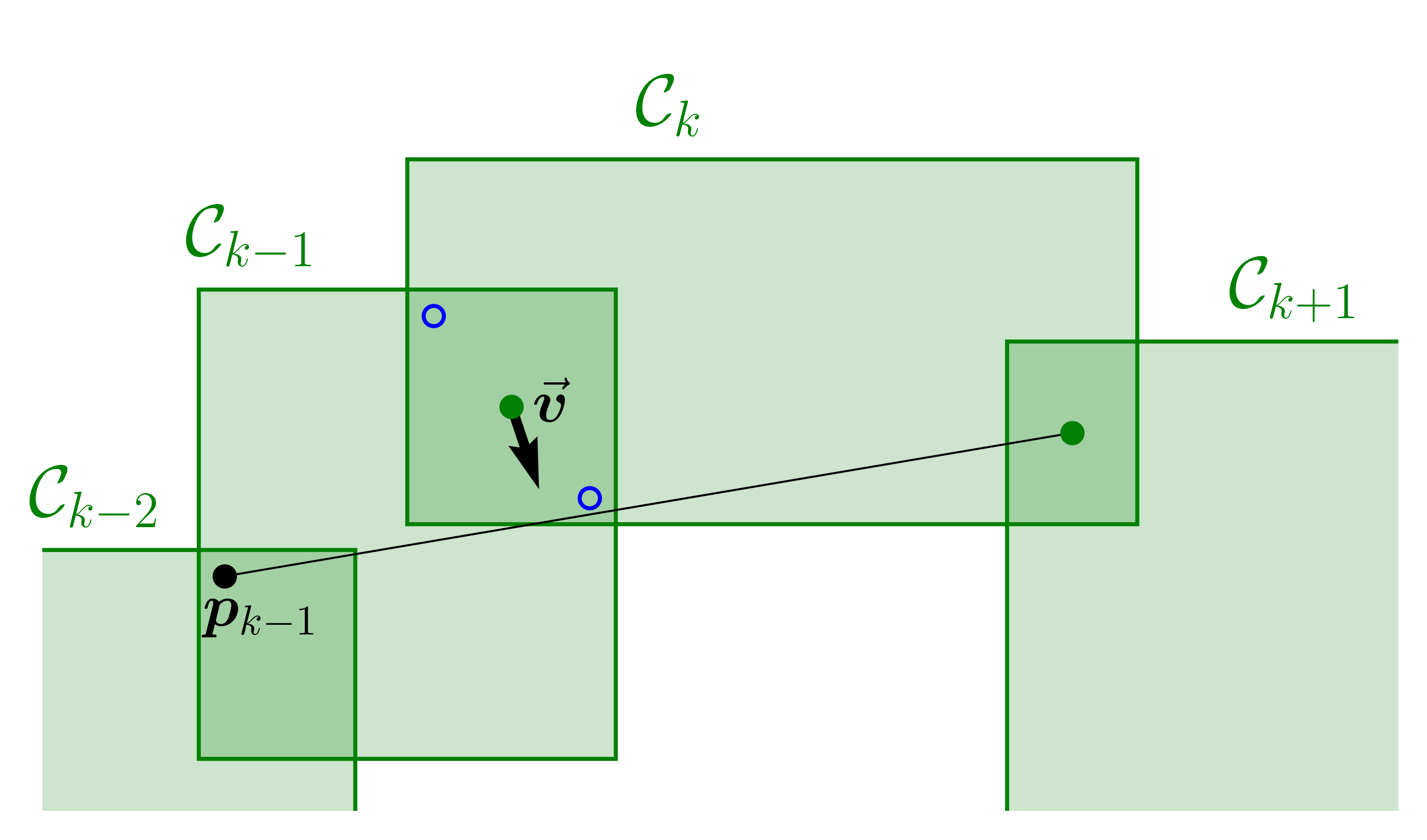}
    \caption{Illustration of heuristics to determine $\{\bm{p}_k\}_{k=1}^{n-1}$. $\bm{C}$ is shown in green and the green dots show centers of $\mathcal{O}_{k-1,k}$ and $\mathcal{O}_{k,k+1}$. The points in the set $\mathcal{W}_k$ are shown as blue circles. The bottom right blue circle will be selected as $\bm{p}_k$}
    \label{fig:waypoint-heuristics}
\end{figure}

The heuristics are based on the assumption that the vehicle takes a turn when moving from $\mathcal{C}_{k-1}$ to $\mathcal{C}_k$. This assumption is fair because the corridors constructed before inflation take 90 degree turns. Therefore, we choose $\{\bm{p}_k\}_{k=1}^{n-1}$ to represent apex points at the inside of the turn. The procedure is illustrated in Figure \ref{fig:waypoint-heuristics}.

The apex point $\bm{p}_k$ is chosen from a set of candidate points denoted as $\mathcal{W}_k$. To construct $\mathcal{W}_k$, we shrink $\mathcal{O}_{k-1,k}$ by $W/2$ on either side horizontally and by $L/2$ on either side vertically. The set $\mathcal{W}_k$ represents the corner points of the shrunken version of $\mathcal{O}_{k-1,k}$. Furthermore, points for which the horizontal distance to an edge of one corridor is larger than $W/2$ and for which the vertical distance to an edge of that same corridor is larger than $L/2$ are removed from $\mathcal{W}_k$ to increase the heuristic's robustness. $\mathcal{W}_k$ is visualized using the blue circles in Figure \ref{fig:waypoint-heuristics}. The heuristic to determine $\bm{p}_k$ states
\begin{equation}
    \bm{p}_k = \argmin_{\bm{q} \in \mathcal{W}_k} || \bm{q} - \bm{i}_k ||_2 \quad k = 1...n-1.
    \label{eq:heuristics-waypoints}
\end{equation}
where $\bm{i}_k$ represents a point at the inside of the turn and is defined by $\bm{i}_k \coloneq \mathrm{center}(\mathcal{O}_{k-1,k}) + \mu \vec{\bm{v}}$ and $\vec{\bm{v}}$ is a unit vector perpendicular to the line connecting $\bm{p}_{k-1}$ and $\mathrm{center}(\mathcal{O}_{k,k+1})$. The sign of $\vec{\bm{v}}$ is chosen such that it points from  $\mathrm{center}(\mathcal{O}_{k-1,k})$ towards this connecting line. The scalar $\mu \in \mathbb{R}_{>0}$ is a dimensionless constant that must be chosen large enough such that $\bm{i}_k \notin \mathcal{O}_{k-1, k}$ to ensure selection of the correct candidate waypoint. As long as this is the case, the exact value of $\mu$ does not affect performance of the method. In this work, we set $\mu = 20$ such that the length of $\mu\vec{\bm{v}}$ has a length of 20m which is larger than the environment dimensions considered in the benchmark and hence larger than any corridor overlap. In words, we select $\bm{p}_k$ as the point in $\mathcal{W}_k$ closest the inside of the turn.

\citet{pontryagin} discuss the analytical solution of the time-optimal trajectory through a series of points $\{\bm{p}_k\}_{k=0}^{n}$ using Pontryagin's Maximum Principle (see Example 4 in their paper). It is shown that the acceleration is a piecewise constant function and that $|a_x(t)| = |a_y(t)| = a_{\mathrm{max}} \quad \forall t \in [0, T]$. The sign of the acceleration is determined by the sign of the co-state trajectory following from the optimality conditions of the Pontryagin's Maximum Principle. However, solving the system of nonlinear equations is nontrivial since it is sensitive to the initialization of the variables of the co-state, which have no direct physical meaning. Additionally, the analysis does not include bounds on the velocity or the position. Therefore, in this work, the sign of the acceleration at points $\{\bm{p}_k\}_{k=0}^{n}$, denoted by $\{\bm{\alpha}_k\}_{k=0}^n$, is heuristically determined as
\begin{equation}
    \bm{\alpha}_k = 
        \begin{cases}
            \mathrm{sign}(\bm{p}_1 - \bm{p}_0) & k=0\\
            \mathrm{sign}(\bm{p}_{n-1} - \bm{p}_f) & k=n\\
            \mathrm{sign}(\bm{p}_k - \mathrm{center}(\mathcal{O}_{k-1,k})) & \mathrm{otherwise.}
        \end{cases}
    \label{eq:heuristics-alpha}
\end{equation}

\begin{figure}
    \centering
    \includegraphics[width=0.8\linewidth]{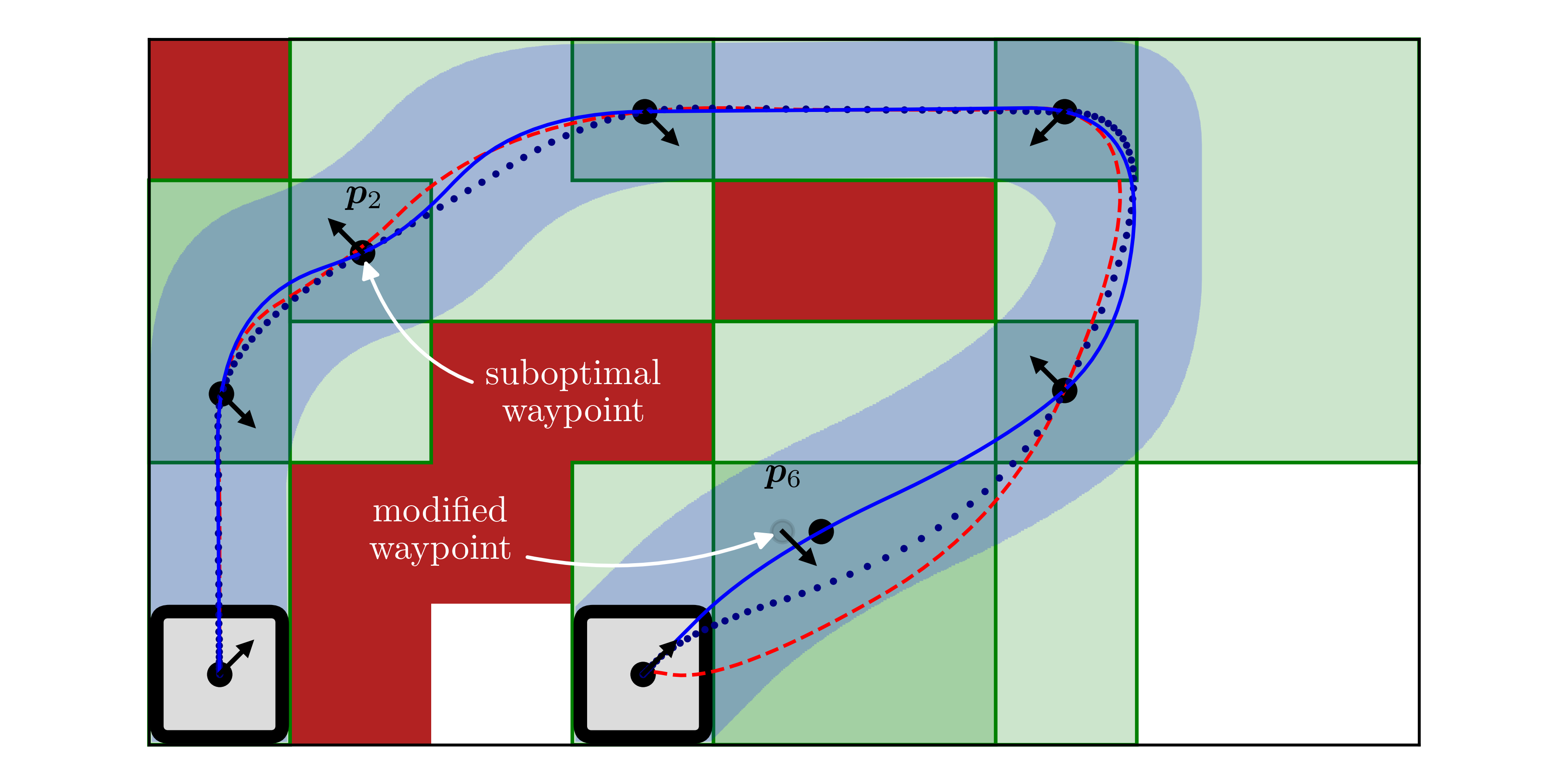}
    \caption{Illustration of selection and modification of $\{\bm{p}_k\}_{k=0}^{n}$ and $\{\bm{\alpha}_k\}_{k=0}^{n}$. $\bm{C}$ is shown in green. The black arrows show $\{\bm{\alpha}_k\}_{k=0}^{n}$. Dashed red line: OCP-solution. Dotted blue line: intermediate PMP-result. Solid blue line: final PMP-solution}
    \label{fig:waypoints-modification}
\end{figure}

In some cases, a modification to $\{\bm{p}_k\}_{k=1}^{n-1}$ and $\{\bm{\alpha}_k\}_{k=1}^{n-1}$ is needed. If the vehicle can move in a straight line from $\bm{p}_{k-1}$ to $\bm{p}_{k+1}$ while staying within the corridors, the assumption that the vehicle makes a turn when moving from one corridor to another is violated. In these cases, the heuristics will produce values $\bm{p}_k$ and $\bm{\alpha}_k$ that are suboptimal or even infeasible. Therefore, the sign of the acceleration is flipped. Formally, if $\bm{q} \in \bigcup_{i=0}^{n} \mathcal{C}_{i}$ $\forall \bm{q} \in \mathcal{F}(\bm{p}_{k-1} + s(\bm{p}_{k+1} - \bm{p}_{k-1}))$ and $\forall s \in [0, 1]$, then $\bm{\alpha}_k \leftarrow -\bm{\alpha}_{k}$ where $\mathcal{F}$ is defined as in \eqref{eq:set-footprint}. Additionally, the point $\bm{p}_k$ is relaxed to be any point in $\mathcal{O}_{k-1,k}$, we set $\bm{p}_k \leftarrow \bm{p}_k + \bm{\Delta}_k$ where $\bm{\Delta}_k$ is an additional constrained degree of freedom, allowing $\bm{p}_k$ to move within the overlapping region. Figure \ref{fig:waypoints-modification} shows an example environment in which a feasible straight line exists from the point $\bm{p}_5$ to the destination, triggering $\bm{\alpha}_6$ to be flipped and $\bm{p}_6$ to be made movable. Note that without this change, $\bm{\Omega}_n = \emptyset$ since no feasible combination of primitives would exist.

\subsection{Parameter Optimization} \label{sec:optimization}
To find the optimal trajectory, $\Omega^*_n = \bm{\Phi}_n(\Pi^*) \in \bm{\Omega}_n$ $(\Pi^* \in \mathcal{X})$ must be found that minimizes the total time of the trajectory
\begin{equation}
    J(\Pi) \coloneq \sum_{k=0}^{n-1}{\left(\sum_{i=0}^{2}{\tau_{i,k}^x}\right)}.
\end{equation} 
To this end, the following optimization problem is solved.

\begin{mini!}|s|
{\bm{z}}{J(\Pi)}
{\label{eq:opti}}{} \label{eq:opti-obj}
\addConstraint{}{\Pi\in \mathcal{X}} \label{eq:opti-feasibility}
\addConstraint{}{0 \leq g_k\left(\bm{p}_k, \bm{v}_k, \bm{\tau}_k\right) \quad}  {k = 0 ... n - 1} \label{eq:opti-corridor_bounds}
\addConstraint{}{|\alpha_0^F| \leq 1} \label{eq:opti-initial-accel}
\addConstraint{}{|\alpha_{n}^F| \leq 1} \label{eq:opti-terminal-accel}
\end{mini!}
where $\bm{z} = \begin{bmatrix}\Pi & \alpha_0^F & \alpha_{n}^F & \bm{\Delta}_1 \dotsi \bm{\Delta}_{n-1}\end{bmatrix}$.

The corridor bounds and velocity limits are enforced using inequality \eqref{eq:opti-corridor_bounds}. To impose velocity limits, the structure of the velocity profile can be exploited. Let $v_k' \coloneq v_k + a_{\mathrm{max}}\alpha_k \tau_{k,0}$ be the velocity obtained after the first acceleration phase. Since the initial velocity profile segment is linear, $v(t) \in [v_k, v_k']$ holds during this first segment. Since the second segment is constant, this inequality still holds. During the third segment (which is again linear), $v(t) \in [v_k', v_{k+1}]$ holds. Hence, it's is sufficient to enforce $\{v_k\}_{k=0}^{n-1} \in [-v_{\mathrm{max}}, v_{\mathrm{max}}]$ and $\{v_k'\}_{k=0}^{n-1} \in [-v_{\mathrm{max}}, v_{\mathrm{max}}]$ leading to $4n$ constraints. Note that $\bm{v}_n$ = 0. The black dots in Figure \ref{fig:primtive-illustration} show these points.

Similarly, the structure of the position profile can be exploited to impose corridor bound constraints sparsely. As long as $|v(t)| > 0$ for the entire duration of the primitive, the position profile is monotonic, meaning feasibility of the primitive is guaranteed if the primitive endpoints are feasible. By construction, this is true since the selected waypoints $\{\bm{p}_k\}_{k=0}^{n}$ all satisfy corridor bounds. If $v(t) = 0$ in one of the parabolic segments, corridor infeasibilities could occur. This is addressed in Section \ref{sec:details}. To prevent the solver from exploring infeasible solutions and hence increasing the robustness, corridor constraints are added to end-points of the coasting stage. This is indicated with the black dots in Figure \ref{fig:primtive-illustration}. Recall that the use of rectangular corridors allow to consider the bounds on position independently in each dimension, which would not be the case if ellipses, such as presented by \citet{large-convex-regions} were used.

\begin{figure}
    \centering
    \includegraphics[width=0.8\linewidth]{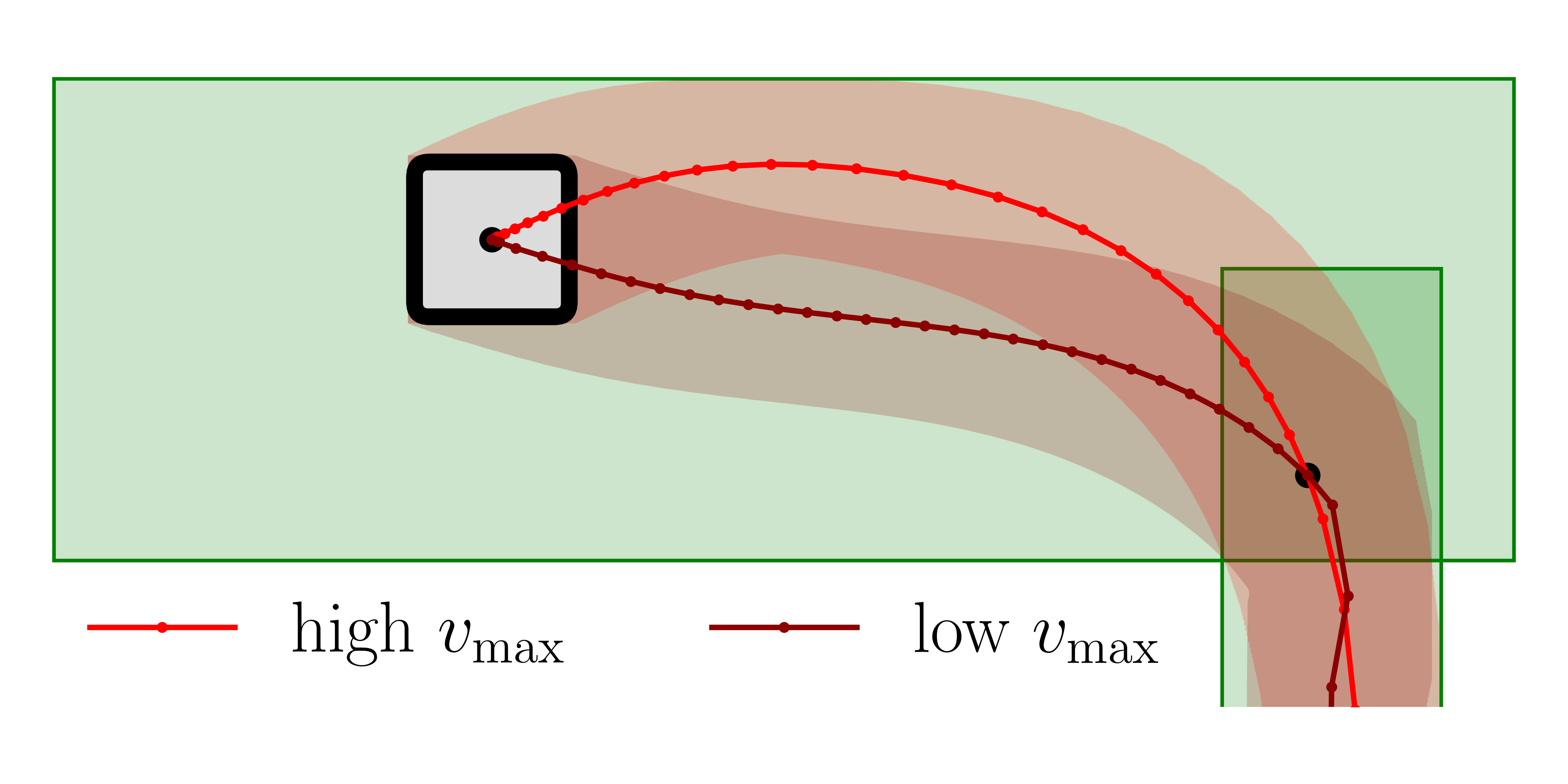}
    \caption{Solution of the OCP for different velocity limits, showing a change in the sign of the initial acceleration in the $y$-direction}
    \label{fig:free-accel}
\end{figure}

Because the sign of the initial acceleration in the direction with the longest distance to travel (referred to as the free direction) is dependent on vehicle parameters, as illustrated in Figure \ref{fig:free-accel}, or even on corridor dimensions, $\alpha^F_0$ is left as an optimization variable. For the same reason, $\alpha^F_n$ is left as an optimization variable. This is captured in \eqref{eq:opti-initial-accel}-\eqref{eq:opti-terminal-accel}.

\subsection{Formulation Details of the Optimization Problem} \label{sec:details}
This section discusses some details in the formulation of the optimization problem that either improve the robustness of the formulation or improve the quality of the solution.

Firstly, the initial and final acceleration in the free direction, respectively denoted by $\alpha_0^F$ and $\alpha_n^F$, are left as optimization variables. Instead of enforcing \eqref{eq:opti-initial-accel} and \eqref{eq:opti-terminal-accel} directly, these constraints are slacked, which was found to help converge to a solution more robustly. The slacked constraints become
\begin{align}
    |\alpha_0^F| &\leq 1 + \left(s_a^0\right)^2, \\
    |\alpha_n^F| &\leq 1 + \left(s_a^n\right)^2.
\end{align}
A penalty $w||\bm{s_a}||_2^2$ where $w = 10^3$ is added to the objective.

Secondly, when problem \eqref{eq:opti} is solved, the resulting trajectory is checked for corridor infeasibilities. If $v(t) = 0$ in one of the parabolic segments, an extreme point in the position profile is reached. If such an extreme point violates the corridor bounds, additional constraints are added and the problem is solved again. For example, suppose the velocity changes sign in the first segment of motion primitive $k$ and $p_k(t)$ violates corridor bounds. In that case, the point $p_k + v_k t + 0.5 \alpha_k a_{\mathrm{max}} t^2$ is constrained to be within corridor bounds as well for $t~=~-v_k / (\alpha_k a_{\mathrm{max}})$. This ensures corridor constraint satisfaction of the extreme points of the position profiles and hence the complete motion primitive. These constraints cannot be added before solving the problem since they would lead to conservatism in many cases.

Thirdly, after obtaining the optimal and feasible trajectory $\Omega^*_n \in \bm{\Omega}_n$, a suboptimal motion primitive can be detected. More precisely, if the vehicle is found coasting in a straight line through a waypoint $\bm{p}_k$, it is not exploiting the acceleration allowed by the primitive, indicating that the imposed acceleration profile does not fit the true time-optimal solution. Therefore, we assign $\bm{\alpha}_k \gets -\bm{\alpha}_k$. Note that the current solution is still feasible for the new profile since it was coasting through $\bm{p}_k$. Problem \eqref{eq:opti} can be solved again for the new motion primitives while the previous solution serves as a feasible initial guess. This process allows to potentially find a more optimal solution and is illustrated in Figure \ref{fig:waypoints-modification}. The first solution is shown as a dotted line and coasts through waypoint $\bm{p}_2$. After changing $\bm{\alpha}_2$, a new and faster trajectory is found.

Finally, the initialization of the primal variables of the optimization problem is nontrivial since there is no intuitive way to initialize the durations of the acceleration. The optimization problem is initialized in every corridor by setting $\tau_{0,k}^x = \tau_{0,k}^y = \tau'_k$, $\tau_{1,k}^x = \tau_{1,k}^y = 7\tau'_k$ and $\tau_{2,k}^x = \tau_{2,k}^y = 0.2\tau_k'$ where $\tau'_k$ is some number larger than 0.05s and $\bm{v}_k$ is chosen such that the primitive reaches the next waypoint $\bm{p}_{k+1}$.

\subsection{Complete Overview} \label{sec:overview}
Algorithm \ref{alg:overview} summarizes PMP. Line \ref{alg-line:corridors} shows the construction of the corridor sequence. Next, in lines \ref{alg-line:analytical-start}-\ref{alg-line:analytical-stop}, the analytical solution is returned if it is feasible. The heuristics are applied to select the primitives in line \ref{alg-line:selection}. The optimization problem \eqref{eq:opti} is then solved in line \ref{alg-line:optimization}. This occurs within a while-loop as explained in Section \ref{sec:details} to ensure corridor feasibility (lines \ref{alg-line:check-solution-start}-\ref{alg-line:check-solution-stop}) and to eliminate suboptimal motion primitives (lines \ref{alg-line:suboptimality-elimination-start}-\ref{alg-line:suboptimality-elimination-stop}).

{\tiny
\begin{algorithm}
    \caption{Overview of PMP}
    \label{alg:overview}
    \small
    \begin{algorithmic}[1]
        \Require $\bm{p}_0$, $\bm{p}_n$, $\bm{\bar{v}}_0$ and a grid map of the environment
        \Ensure Time-Optimal, collision-free trajectory

        \State construct corridor sequence $\bm{C}$ \label{alg-line:corridors} \Comment{Sec. \ref{sec:corridors}}
        \If {analytical solution is feasible}  \label{alg-line:analytical-start} \Comment{Sec. \ref{sec:analytical}}
            \State \Return analytical solution 
        \EndIf                                 \label{alg-line:analytical-stop}
        \State determine $\{\bm{p}_k\}_{k=1}^{n-1}$ and $\{\bm{\alpha}_k\}_{k=0}^{n}$ \label{alg-line:selection} \Comment{Sec. \ref{sec:primitive-selection}}
        \State completed $\gets$ False
        \While{not completed}
            \State solve problem \eqref{eq:opti}        \label{alg-line:optimization}\Comment{Sec. \ref{sec:optimization}}
            \While{corridor bounds violated}            \label{alg-line:check-solution-start}\Comment{Sec. \ref{sec:details}}
                \State add additional constraints
                \State solve problem \eqref{eq:opti}
            \EndWhile                                   \label{alg-line:check-solution-stop}
            \If{suboptimal primitive found}          \label{alg-line:suboptimality-elimination-start}\Comment{Sec. \ref{sec:details}}
                \State modify primitive selection
            \Else
                \State completed $\gets$ True
            \EndIf                                      \label{alg-line:suboptimality-elimination-stop}
        \EndWhile
        \State \Return solution
    \end{algorithmic}
\end{algorithm}
}

\subsection{Method Limitations} \label{sec:limitations}
The proposed method has some limitations that are addressed in this section. Firstly, the corridor sequence is based on a distance-optimal graph search, which can lead to corridor sequences that are not time-optimal as shown in Figure \ref{fig:unstructured-envs-a} where VP-STO finds a better trajectory outside of the corridor sequence. One way to address this limitation is to explore initial paths $\mathcal{P}$ that belong to different homotopy classes. These different classes can be found by constructing a Voronoi diagram on top of the map of the environment, such as shown by \citet{voronoi}.

Secondly, it is assumed that the grid cells in the environment representation are large enough to contain the vehicle. This ensures that the vehicle can always transition from corridor to corridor but also puts a lower limit on obstacle sizes. Since the proposed method mainly targets structured environments, this assumption is considered fair. Note that for the specific application considered, obstacles will usually be environment edges, faulty tiles or stationary movers, which all fit the assumption. However, in randomly cluttered environments with irregular obstacles, sufficient overlap between corridors cannot be achieved by making this assumption and hence, must be achieved in another way. Additionally, the use of rectangular corridors might not be justified in those environments and might not be able to properly capture the available free space.

Thirdly, the heuristics to select the motion primitives do not provide theoretical guarantees. Indeed, as shown in Figure \ref{fig:unstructured-envs-b}, the heuristics can fail and produce suboptimal trajectories. Additionally, the heuristics do not guarantee that a feasible set of parameters exists. However, benchmarking results in Section \ref{sec:benchmark} demonstrate that the heuristics are effective in practice.

Finally, equation \eqref{eq:velocity-acceleration-limits} imposes velocity and acceleration bounds using the inifinity-norm. This is accurate if the motion of the vehicle is fully independent in each dimension, but can introduce suboptimality if the total velocity or acceleration (Euclidian norm) is limited rather than the infinity norm. Even still, the choice of this norm allows to decouple acceleration profiles in each dimension, which in turn allows to significantly reduce computation times.
\section{Benchmarking} \label{sec:benchmark}
To evaluate the proposed method (PMP), we compare against the methods described in Section \ref{sec:benchmark-method-overview}. All methods are applied to compute trajectories in a structured environment and in randomly cluttered environments. Notes on implementation details are provided in Section \ref{sec:implementation} and the benchmark itself is detailed in Section \ref{sec:benchmark-description}. Finally, we compare the different methods to construct dynamically feasible paths in Section \ref{subsec:benchmark-optimization}.

\subsection{Overview of Methods Benchmarked Against} \label{sec:benchmark-method-overview}
The proposed method is compared against two other optimization-based approaches that plan through the same corridor sequence $\bm{C}$. This allows to validate the selection and optimization of the parametric motion primitives within the corridors. In addition, we compare against VP-STO which does not use the corridor sequence $\bm{C}$.

\subsubsection{Full Optimal Control Problem (OCP)}
A full optimization problem can be written down in which the trajectory is discretized using a multiple-shooting approach. The Optimal Control Problem (OCP) is defined using a multi-stage formulation with time-scaling where every stage represents the movement through one corridor. The velocity and position bounds are applied to every discretization point. The formulation is similar to problem (3) in the work of \citet{truck-trailer}. For fair comparison, we use an analytical integrator for the dynamics. The number of grid-points per corridor affects the size of the problem, as well as the feasibility and optimality of the trajectory. A finer grid is better suited to capture the bang-bang solution and is less likely to produce a trajectory that violates corridor bounds in between grid-points. However, it comes at the cost of resulting in a larger problem with higher solver times. In this work, the use of 30 grid-points per corridor was found to be a fair trade-off.

\subsubsection{OMG-tools}
OMG-tools is a toolbox that uses spline representations to produce smooth trajectories \citep{omg2, omg1}. Both the position and the velocity are represented using a spline and constraints on the acceleration are applied to the derivative of the velocity.

\subsubsection{Via-point-based Stochastic Trajectory Optimization (VP-STO)} \label{sec:vp-sto-intro}
Instead of planning through the corridor sequence, VP-STO \citep{vp-sto} samples a number of points ($N_{\mathrm{via}}$) in space and connects these with the analytical solution of an OCP. Using stochastic black-box and gradient-free optimization, the via-points are modified to minimize some general cost-function. By increasing ($N_{\mathrm{via}}$), more complex trajectories can be represented. Additionally, by increasing the population size, meaning the number of trajectories considered, the probability of obtaining the globally optimal solution is increased. However, by increasing the value of these parameters, so does the computation time. In this benchmark, we consider a population size of 25, max iteration count of 3000, $\sigma_{\mathrm{init}}$ equal to 1.5 and different values of $N_\mathrm{via}$. We mention this number by referring to this method as VP-STO-$N_\mathrm{via}$.

\citet{vp-sto-git} provide an open-source implementation in Python. However, this implementation is only intended for offline use and can likely be significantly sped up. The authors also propose online-VP-STO which they use for real-world experiments. In that setting, only a limited number of via-points is considered and the trajectory is not solved until convergence but is instead aborted when the MPC-deadline is reached, sacrificing optimality in the process. Additionally, the authors state that when the target location was changed, the robot was sometimes unable to find a valid motion quickly enough, indicating a reliance on the warm-starting procedure. In this work, we consider VP-STO and run it until convergence (or until the maximum iteration count is reached).

\subsection{Implementation Specifics} \label{sec:implementation}
All optimization problems are formulated using \texttt{CasADi} \citep{casadi}. For OCP, an opti function object is created offline for every amount of corridors to be considered. These are solved using FATROP \citep{vanroye2023fatrop}, a state-of-the-art OCP-solver based on IPOPT \citep{ipopt}.

Similarly, opti functions are created for PMP. However, for $n$ corridors, there are $n-1$ possible movable waypoints which requires enumerating $2^{n-1}$ variations of the optimization problem. However, these problems are small and there is no need to enumerate every possible variation as long as a variation that allows for enough movable points is available. These optimization problems are also solved using FATROP.

OMG-tools is only available in \texttt{Python} and we solve the optimization problems using IPOPT with linear solver MA27. The code is executed on a 12th Gen Intel(R) Core(TM) i7-1280P processor. Note that the implementation of OMG-tools can likely still be sped up. However, the optimization problem formulated by OMG-tools currently does not fit the problem structure required by FATROP and hence must be solved using IPOPT.

\subsection{Benchmark Description} \label{sec:benchmark-description}
The methods considered are applied to structured and randomly cluttered environments. Figure \ref{fig:structured-envs} shows the structured environment. It contains some aligned walls and obstacles. Within this environment, 100 trajectories with different parameter values, starting positions and destinations are computed. The parameters $v_{\mathrm{max}}$ and $a_{\mathrm{max}}$ are sampled uniformly from $[0.5, 2.0] m/s$ and $[2.0, 6.0] m/s^2$ respectively. Feasible points $\bm{p}_0$ and $\bm{p}_n$ are selected uniformly as well but are only accepted if $||\bm{p}_0 - \bm{p}_n||_2 > 5W$ where $W$ is the vehicle width (which is equal in this case to its length) to avoid unreasonably short trajectories.

Apart from this structured environment, 100 random environments are created in which the probability of a grid cell containing and obstacle is set to 0.1. For higher obstacle densities, similar results were found. Parameter values, starting positions and terminal positions are randomly selected as described for the structured environment.

\subsection{Results} \label{subsec:benchmark-optimization}
\begin{table*}
    \centering
    \caption{Comparison of 100 trajectories with random starting positions, destinations and values for $a_\mathrm{max}$ and $v_\mathrm{max}$}
    \label{tab:my_label}
    \begin{subtable}{\linewidth}
        \subcaption{Structured environment}
        \label{tab:structured}
        \begin{tabular}{r|ccc|ccc}
            \toprule
                 & PMP (ours) & OCP & OMG-tools & VP-STO-5 & VP-STO-10 & VP-STO-15 \\
                \midrule
                avg $t_\mathrm{solver}$ [ms] & \textbf{2.38} & 11.33 & 12.03 & (8.28e+03) & (2.30e+04) & (4.92e+04) \\
                max $t_\mathrm{solver}$ [ms] & \textbf{6.61} & 28.09 & 21.58 & (3.12e+04) & (1.22e+05) & (1.57e+05) \\
                avg $t_\mathrm{total}$ [ms] & \textbf{4.30} & 15.09 & 146.15 & (8.42e+03) & (2.31e+04) & (4.94e+04) \\
                max $t_\mathrm{total}$ [ms] & \textbf{9.83} & 33.18 & 710.91 & (3.17e+04) & (1.22e+05) & (1.57e+05) \\ [1em]
                avg $t_\mathrm{move}$ [s] & 1.844 & \textbf{1.839} & 1.957 & 2.086 & 1.894 & 1.858 \\
                $\varepsilon_\mathrm{move}$ [\%] & $0.0 \pm 1.3$ & - & $6.0 \pm 4.9$ & $6.7 \pm 14.6$ & $1.9 \pm 11.8$ & $0.7 \pm 11.3$ \\ [1em]
                avg $t_\mathrm{total} + t_\mathrm{move}$ [s] & \textbf{1.849} & 1.854 & 2.103 & 10.508 & 25.006 & 51.249 \\ [1em]
                \# infeasible cases & \textbf{0} & 13 & \textbf{0} & \textbf{0} & \textbf{0} & \textbf{0} \\
                \# solver failures & 1 & \textbf{0} & 1 & 5 & 6 & 10 \\
                \bottomrule
            \end{tabular}
    \end{subtable}

    \begin{subtable}{\linewidth}
        \subcaption{Unstructured environment}
        \label{tab:unstructured}
        \begin{tabular}{r|ccc|ccc}
            \toprule
                 & PMP (ours) & OCP & OMG-tools & VP-STO-5 & VP-STO-10 & VP-STO-15 \\
                \midrule
                avg $t_\mathrm{solver}$ [ms] & \textbf{4.28} & 16.68 & 16.89 & (3.13e+03) & (8.93e+03) & (1.85e+04) \\
                max $t_\mathrm{solver}$ [ms] & \textbf{26.10} & 60.29 & 36.43 & (2.27e+04) & (4.01e+04) & (5.53e+04) \\
                avg $t_\mathrm{total}$ [ms] & \textbf{6.42} & 20.90 & 181.51 & (3.25e+03) & (9.05e+03) & (1.86e+04) \\
                max $t_\mathrm{total}$ [ms] & \textbf{29.00} & 64.37 & 648.17 & (2.29e+04) & (4.02e+04) & (5.54e+04) \\ [1em]
                avg $t_\mathrm{move}$ [s] & 1.709 & \textbf{1.703} & 1.785 & 1.893 & 1.724 & 1.694 \\
                $\varepsilon_\mathrm{move}$ [\%] & $0.0 \pm 1.7$ & - & $4.2 \pm 3.2$ & $6.1 \pm 19.9$ & $0.9 \pm 13.0$ & $0.5 \pm 11.2$ \\ [1em]
                avg $t_\mathrm{total} + t_\mathrm{move}$ [s] & \textbf{1.716} & 1.724 & 1.966 & 5.146 & 10.774 & 20.297 \\ [1em]
                \# infeasible cases & \textbf{0} & 7 & \textbf{0} & \textbf{0} & \textbf{0} & \textbf{0} \\
                \# solver failures & 3 & \textbf{0} & 4 & 5 & 3 & 3 \\
                \bottomrule
            \end{tabular}
        
    \end{subtable}
\end{table*}

For every trajectory in the benchmark, we consider the time to compute the trajectory, the moving time of the computed trajectory (which is minimized), and the robustness of the method. The results are summarized in Table \ref{tab:my_label} and explained in detail in the next subsections.

\subsubsection{Discussion of Computation Time}
\begin{figure}
    \centering
    \begin{subfigure}{0.9\linewidth}
        \includegraphics[width=\linewidth]{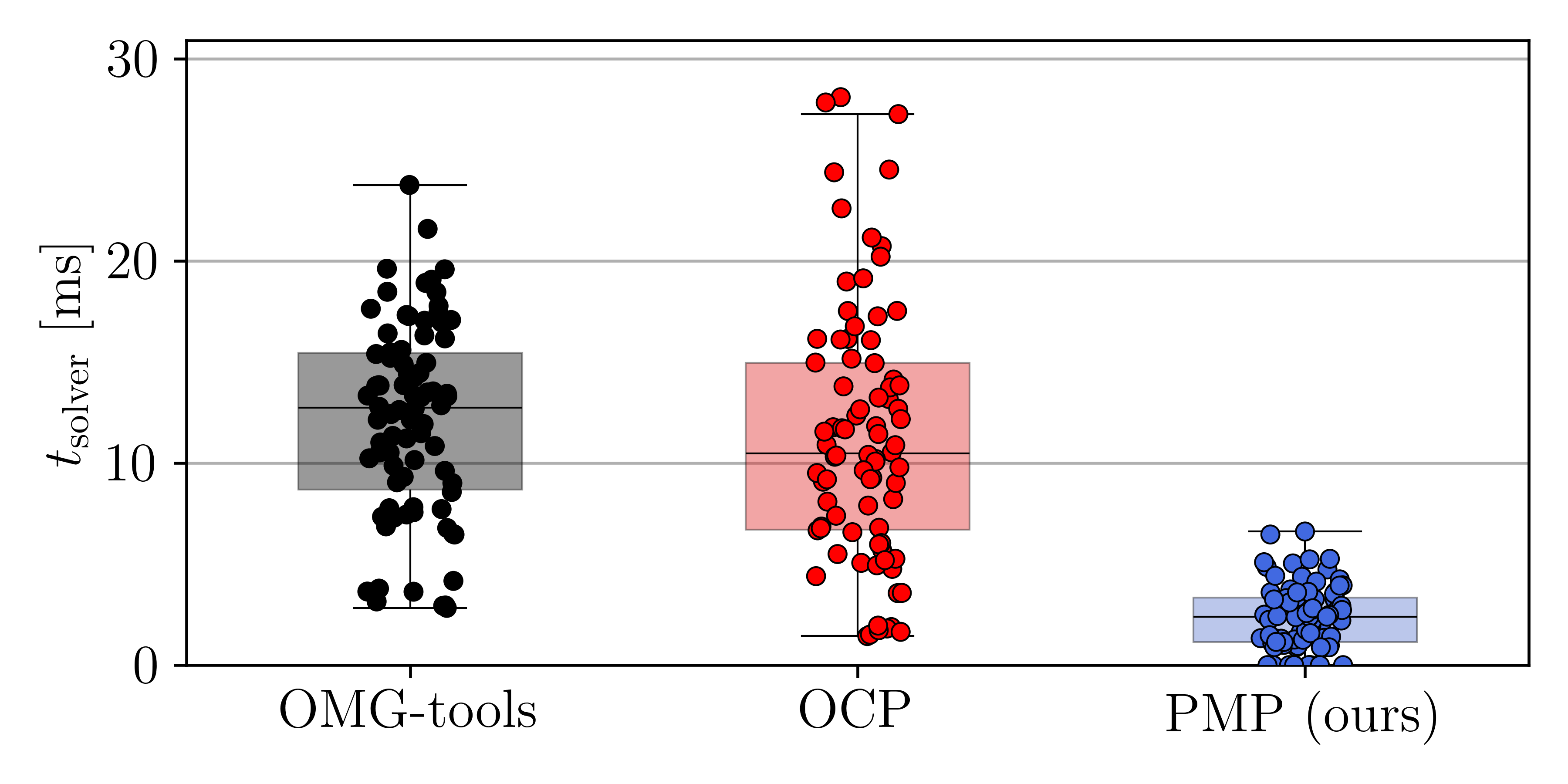}
        \subcaption{Structured environment}
        \label{fig:solver-times-structured}
    \end{subfigure}
    \begin{subfigure}{0.9\linewidth}
        \includegraphics[width=\linewidth]{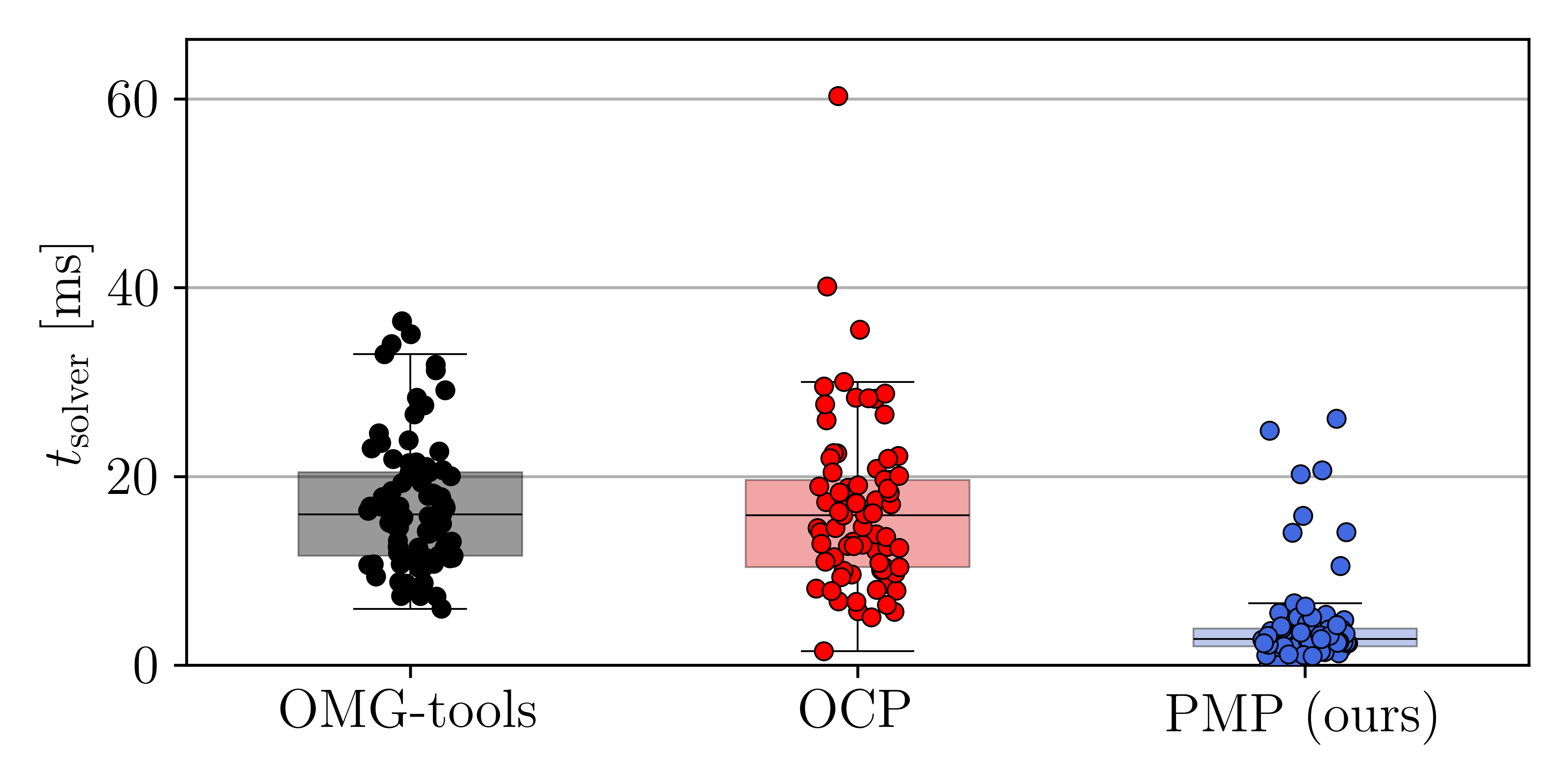}
        \subcaption{Unstructured environment}
        \label{fig:solver-times-unstructured}
    \end{subfigure}
    \caption{Comparison of solver times ($t_{\mathrm{solver}}$) for the optimization-based methods on the benchmark}
    \label{fig:solver-times}
\end{figure}

The total computation time ($t_\mathrm{total}$) includes the solver time ($t_{\mathrm{solver}}$), which is the time needed to solve the optimization problem but also includes corridor construction, potentially the selection of motion primitives, setting initial guesses and sampling the trajectory. Table \ref{tab:my_label} shows average and worst case values for $t_{\mathrm{solver}}$ and $t_\mathrm{total}$ for both type of environments. PMP is able to solve the optimization problem about 4 times faster than OCP or OMG-tools. Note that OMG-tools has a significantly higher total time, which is due to inefficient construction of the optimization problem, and could be improved with a better implementation. Figure \ref{fig:solver-times} shows the solver times for the corridor-based methods. In the unstructured environment, high solver times are observed for PMP in some cases. This is caused by the fact that the implicit assumption that the vehicle turns when transitioning from one corridor to the next is not always valid in the unstructured environment. This can cause a suboptimal selection of primitives which is detected and corrected by solving the problem again, at the cost of a higher computation time.

In the structured environment, PMP was able to solve the problem using the analytical solution to plan through the corridor sequence described in Section \ref{sec:analytical} in 10 cases. In 8 out of those 10 cases, the corridor sequence consisted of a single corridor. In the remaining 2, two corridors made up the corridor sequence. In the unstructured environment, PMP used the analytical solution in 5 cases. In 4 of those cases the corridor sequence consisted of two corridors and in only 1 case there was a single corridor.

The computation times for VP-STO are significantly higher than those of the corridor-based methods. In addition, more via-points lead to higher computation times. These results are not fully representative since the implementation can likely be sped up significantly as mentioned in Section \ref{sec:vp-sto-intro}.

\subsubsection{Discussion of Optimality}
\begin{figure}
    \centering
    \begin{subfigure}{0.9\linewidth}
        \includegraphics[width=\linewidth]{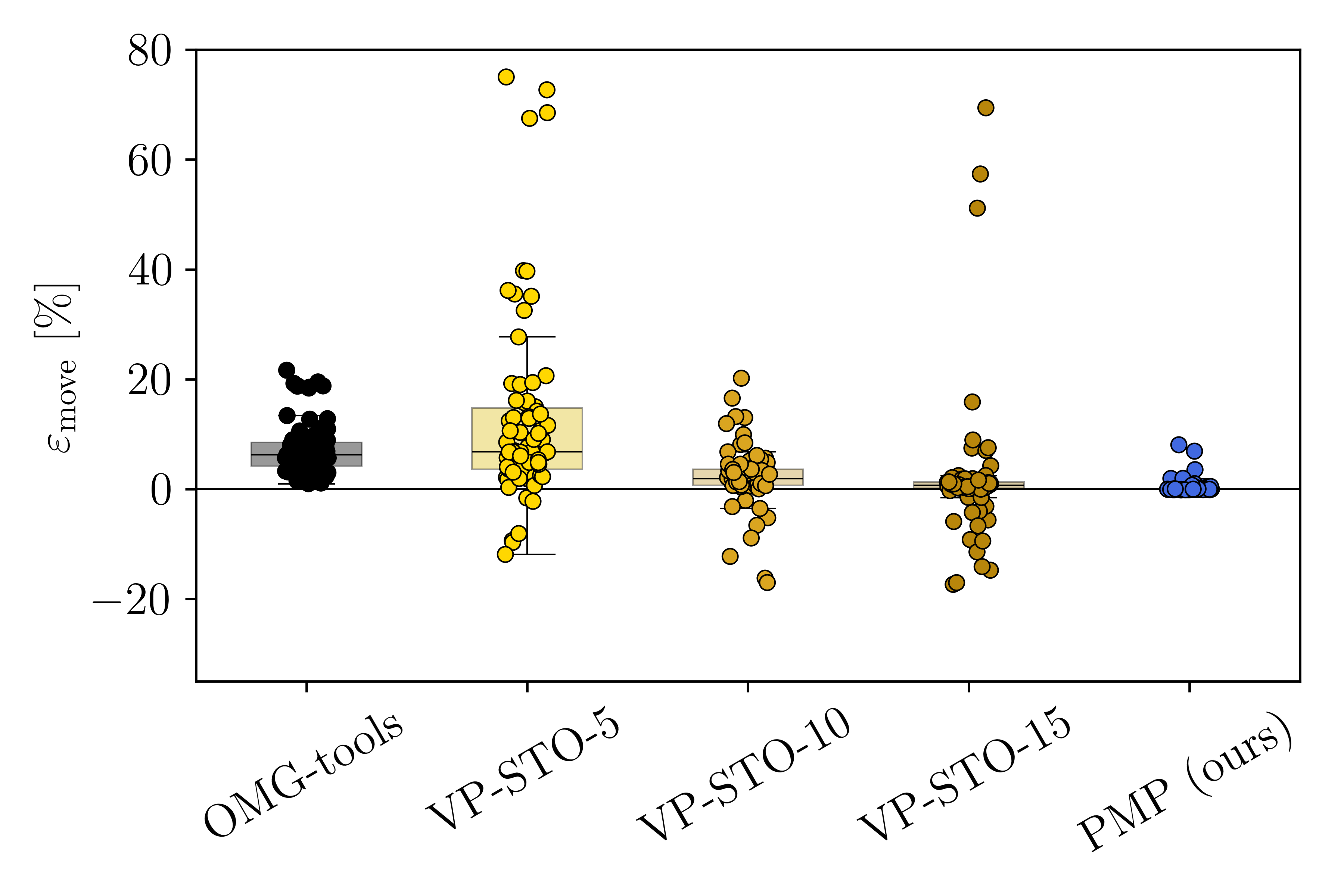}
        \subcaption{Structured environment}
        \label{fig:optimality-structured}
    \end{subfigure}
    \begin{subfigure}{0.9\linewidth}
        \includegraphics[width=\linewidth]{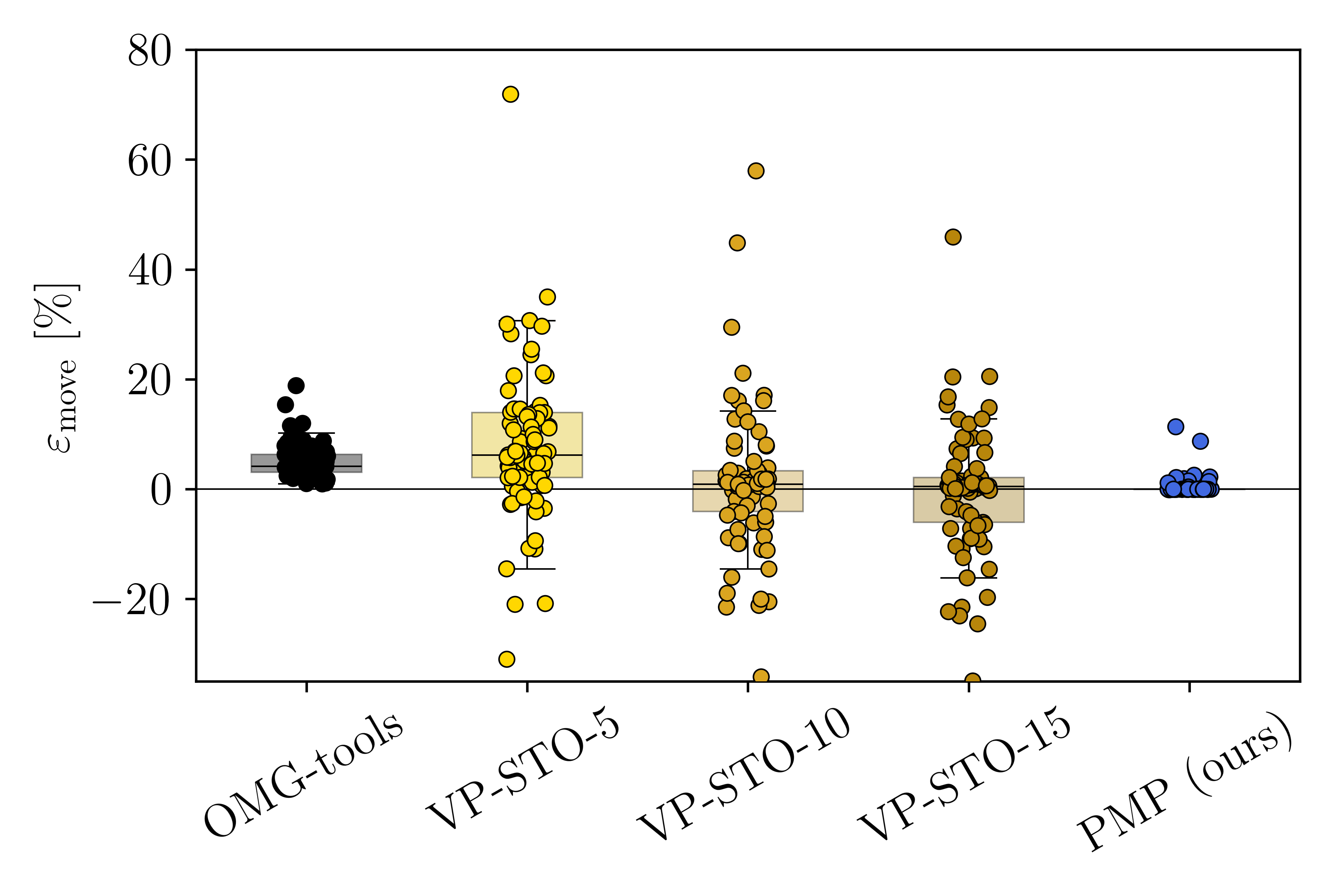}
        \subcaption{Unstructured environment}
        \label{fig:optimality-unstructured}
    \end{subfigure}
    \caption{Relative error on the travel time for all environments. OCP is considered the baseline}
    \label{fig:suboptimality}
\end{figure}

\begin{figure}
    \centering
    \begin{subfigure}{0.45\linewidth}
        \includegraphics[width=\linewidth]{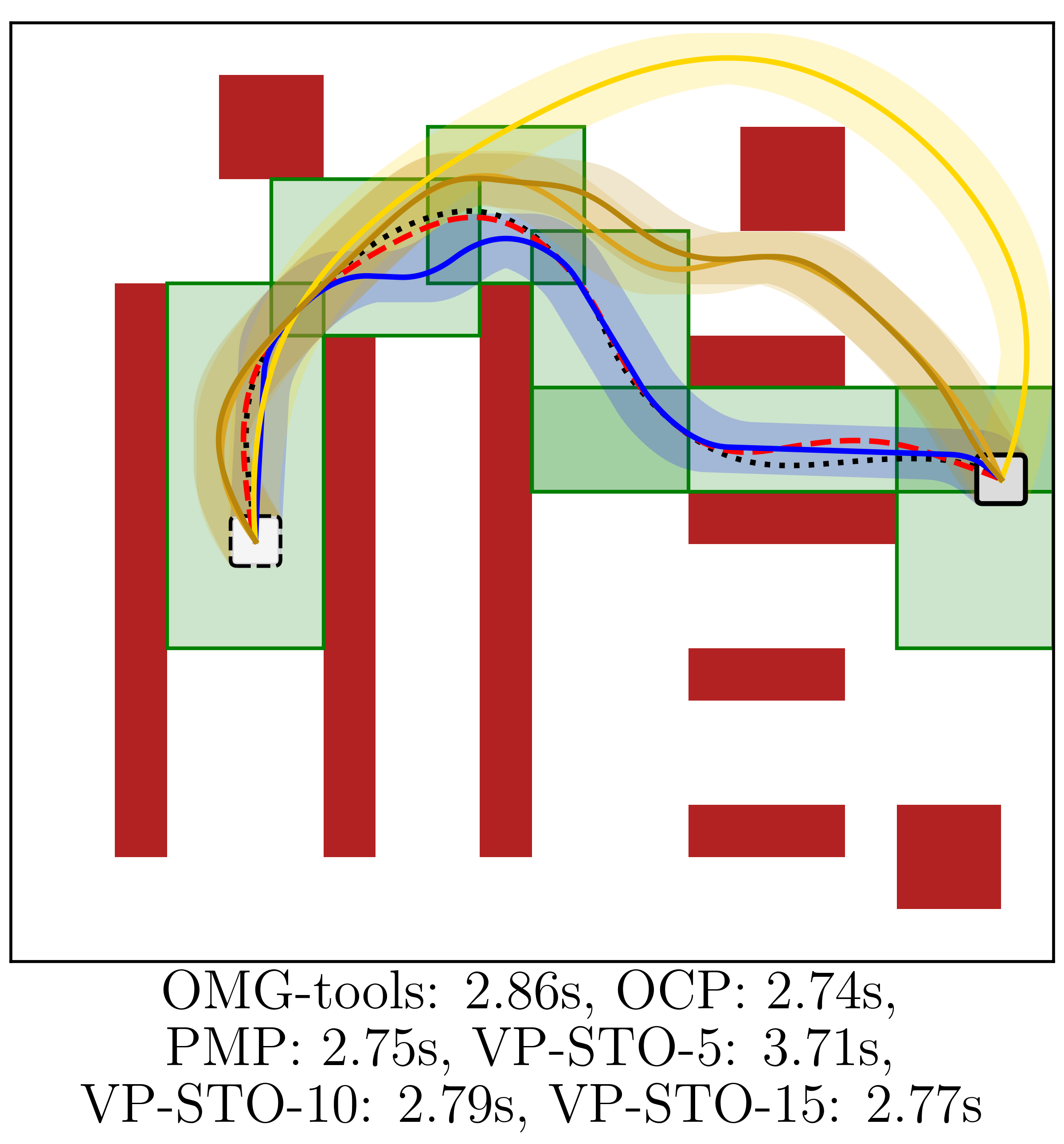}
        \subcaption{}
        \label{fig:structured-envs-a}
    \end{subfigure}
    \begin{subfigure}{0.45\linewidth}
        \includegraphics[width=\linewidth]{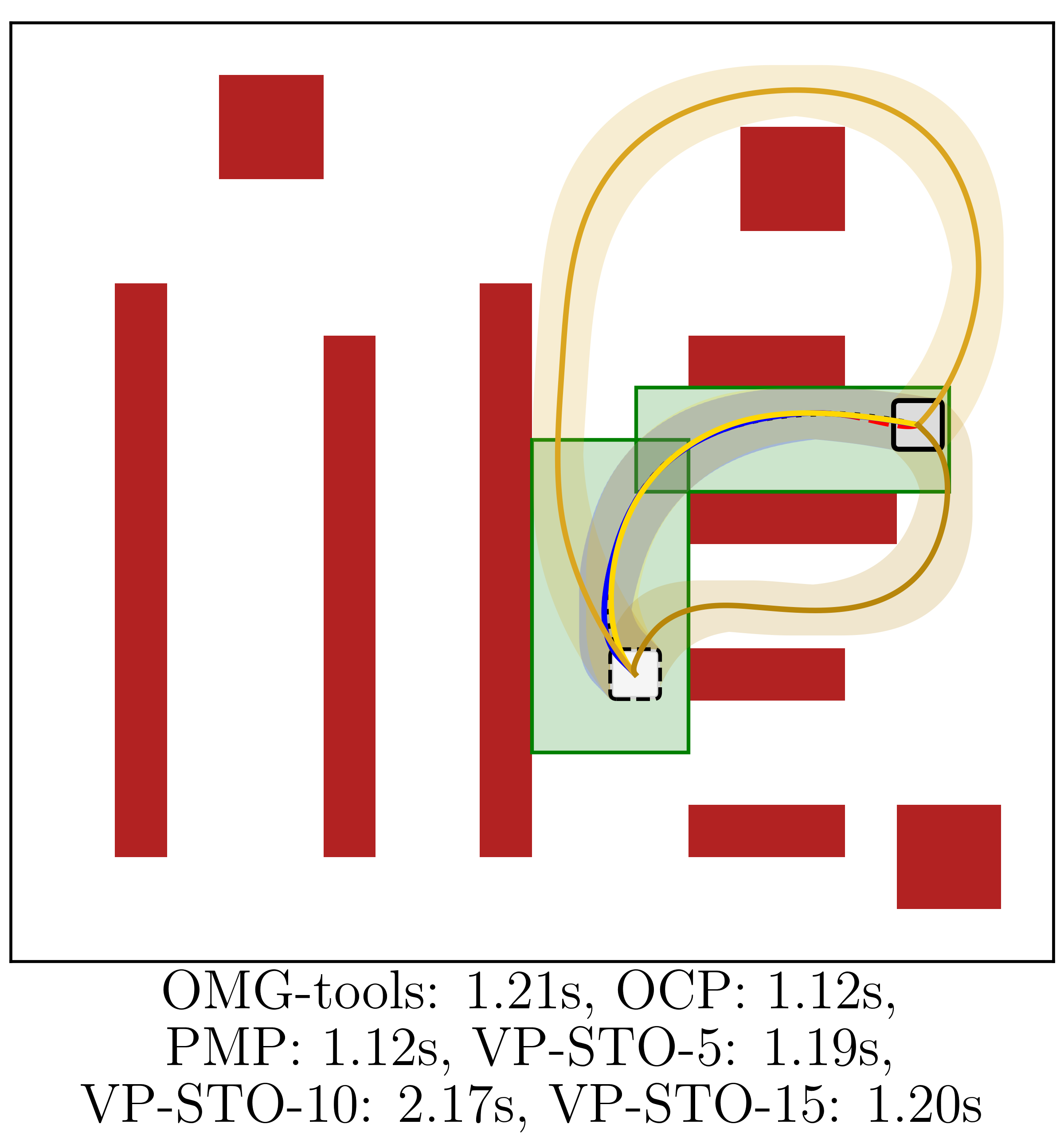}
        \subcaption{}
        \label{fig:structured-envs-b}
    \end{subfigure}
    \begin{subfigure}{0.45\linewidth}
        \includegraphics[width=\linewidth]{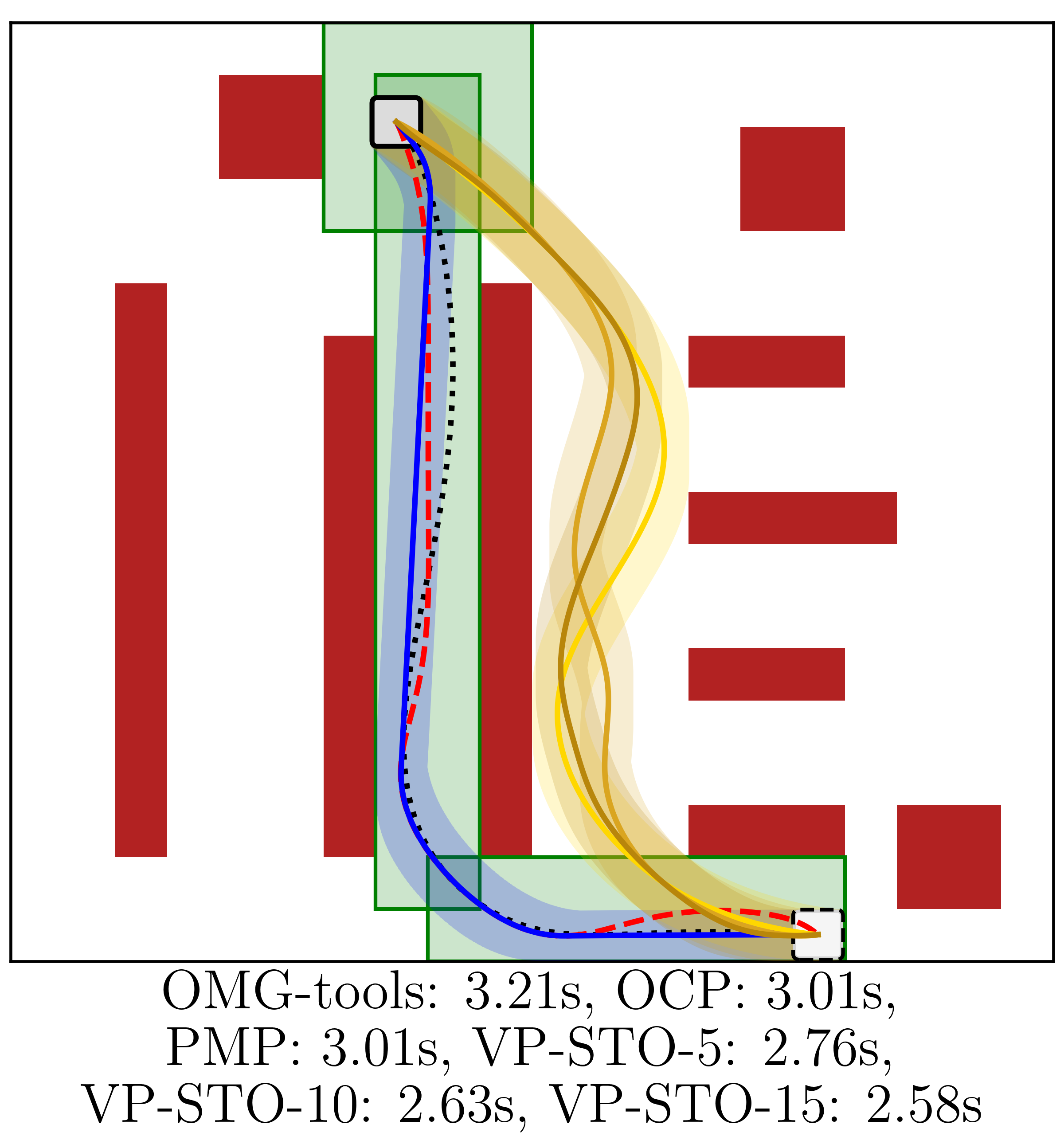}
        \subcaption{}
        \label{fig:structured-envs-c}
    \end{subfigure}
    \begin{subfigure}{0.45\linewidth}
        \includegraphics[width=\linewidth]{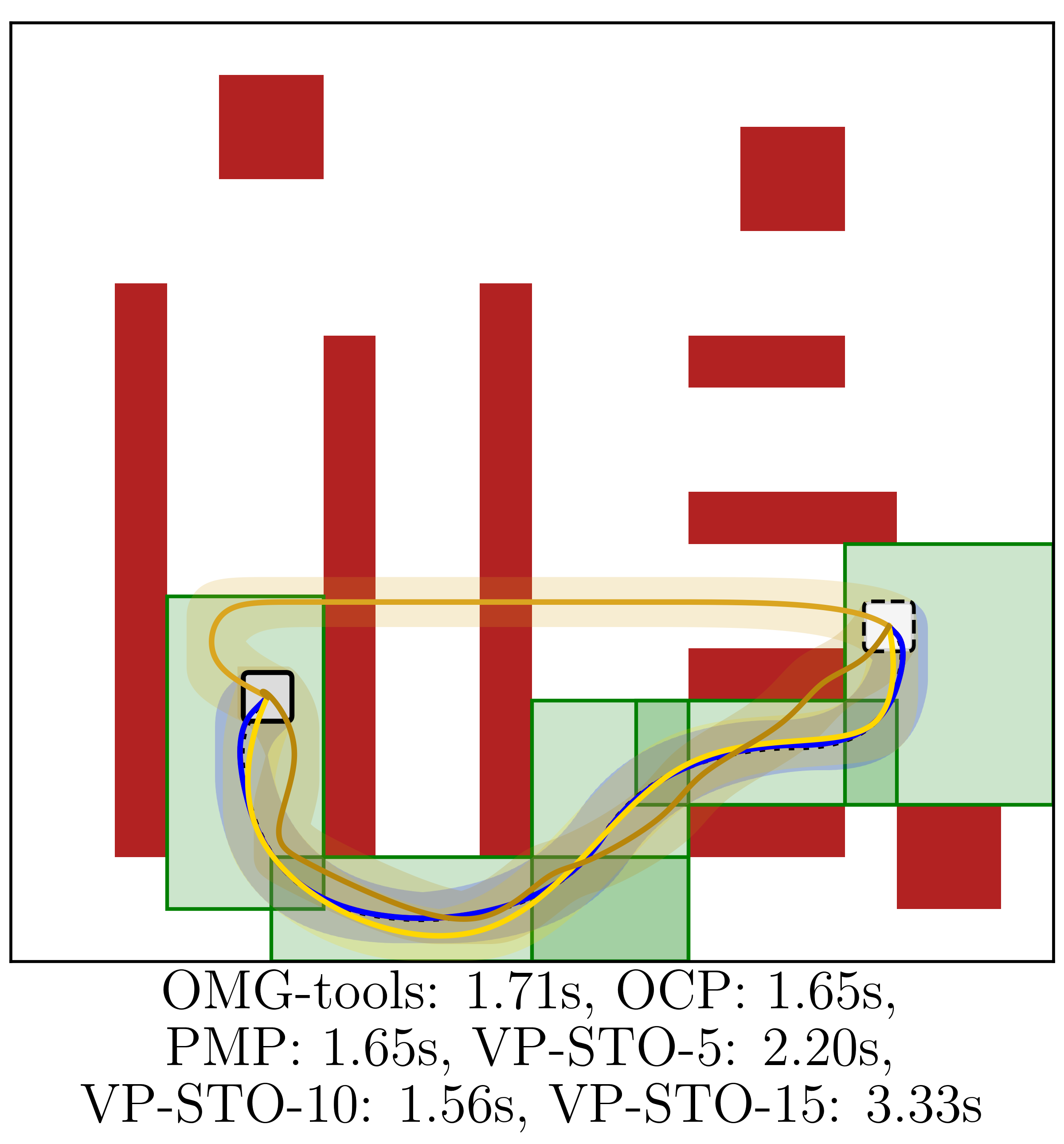}
        \subcaption{}
        \label{fig:structured-envs-d}
    \end{subfigure}
    \begin{subfigure}{\linewidth}
        \includegraphics[width=\linewidth]{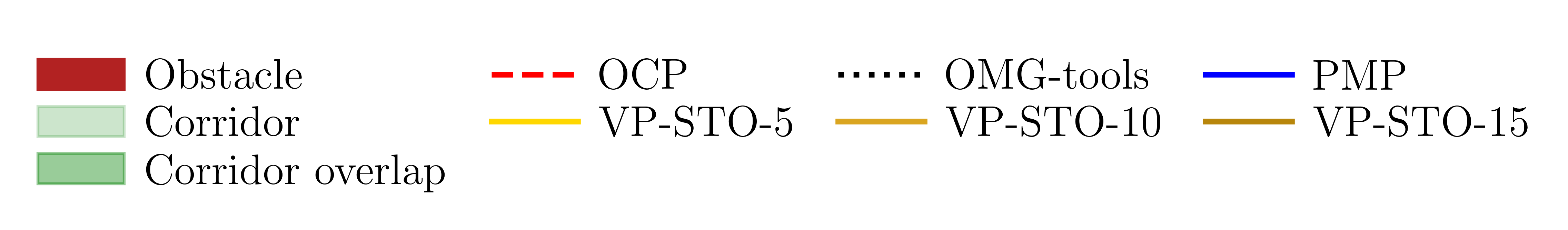}
    \end{subfigure}
    \caption{Some example trajectories in the structured environment. The vehicle's footprint at $\bm{p}_0$ (solid) and $\bm{p}_n$ (dashed) is shown. Transparently shaded regions show traces of the vehicle footprint}
    \label{fig:structured-envs}
\end{figure}

\begin{figure}
    \centering
    \begin{subfigure}{0.45\linewidth}
        \includegraphics[width=\linewidth]{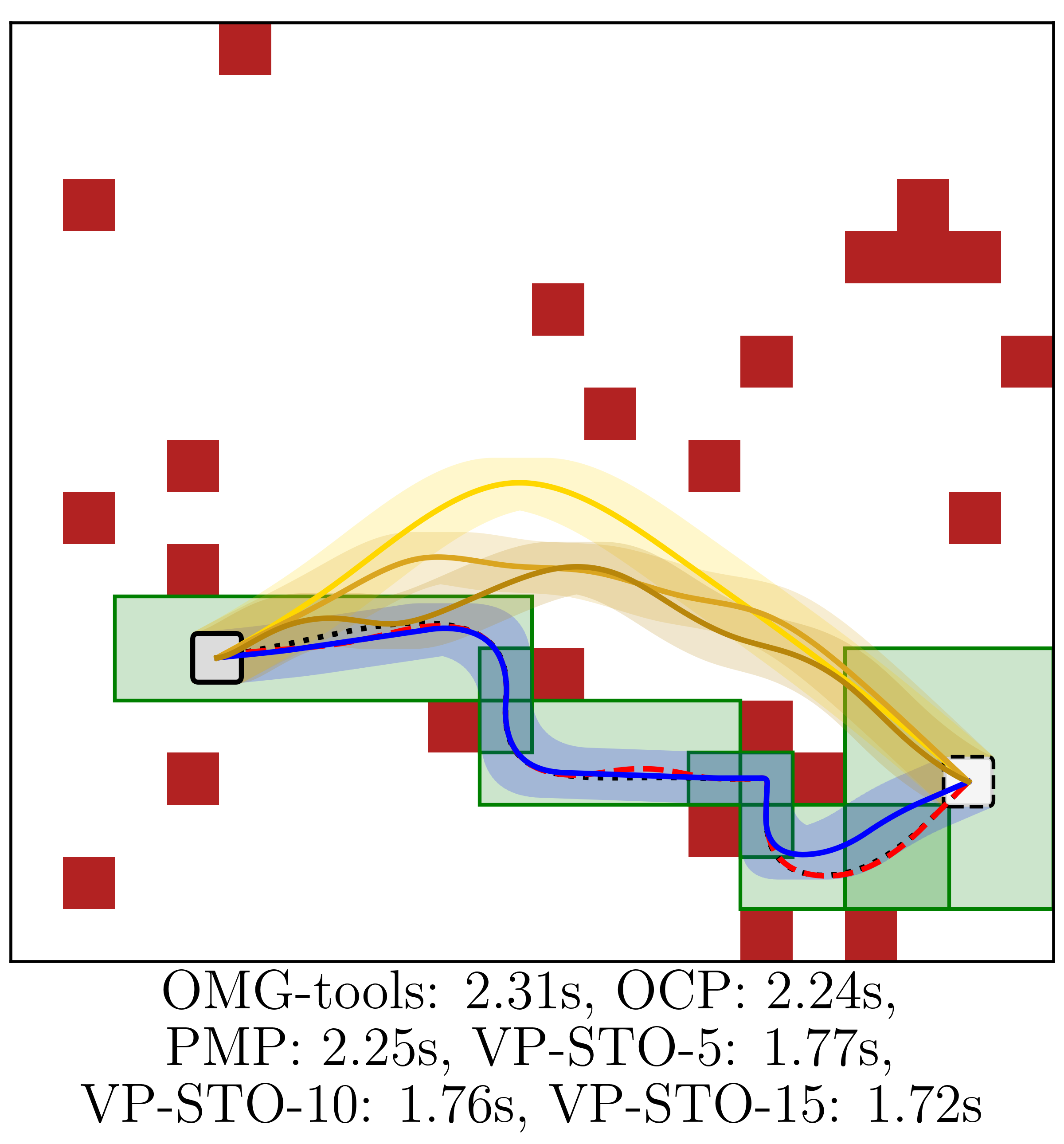}
        \subcaption{}
        \label{fig:unstructured-envs-a}
    \end{subfigure}
    \begin{subfigure}{0.45\linewidth}
        \includegraphics[width=\linewidth]{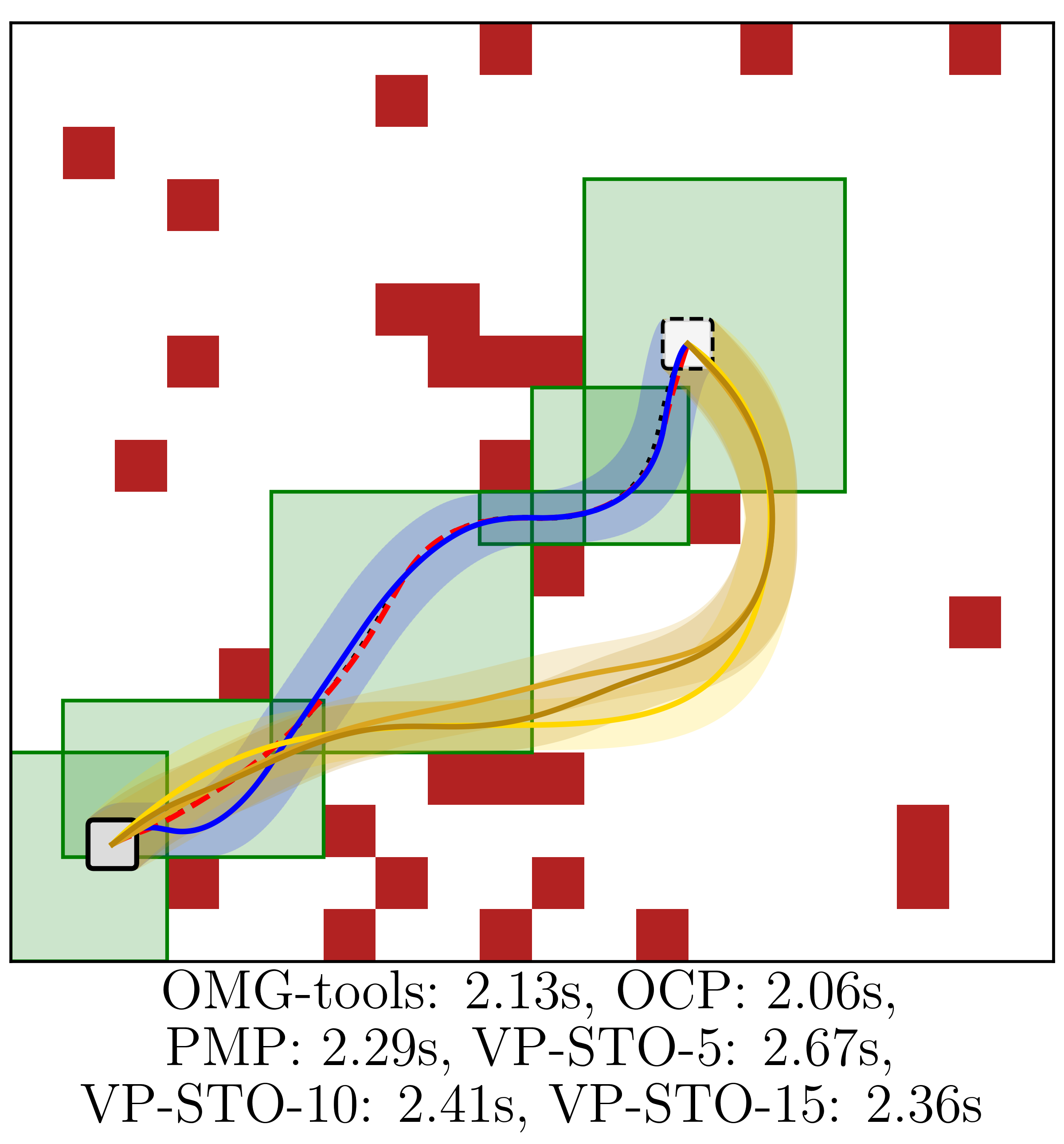}
        \subcaption{}
        \label{fig:unstructured-envs-b}
    \end{subfigure}
    \begin{subfigure}{0.45\linewidth}
        \includegraphics[width=\linewidth]{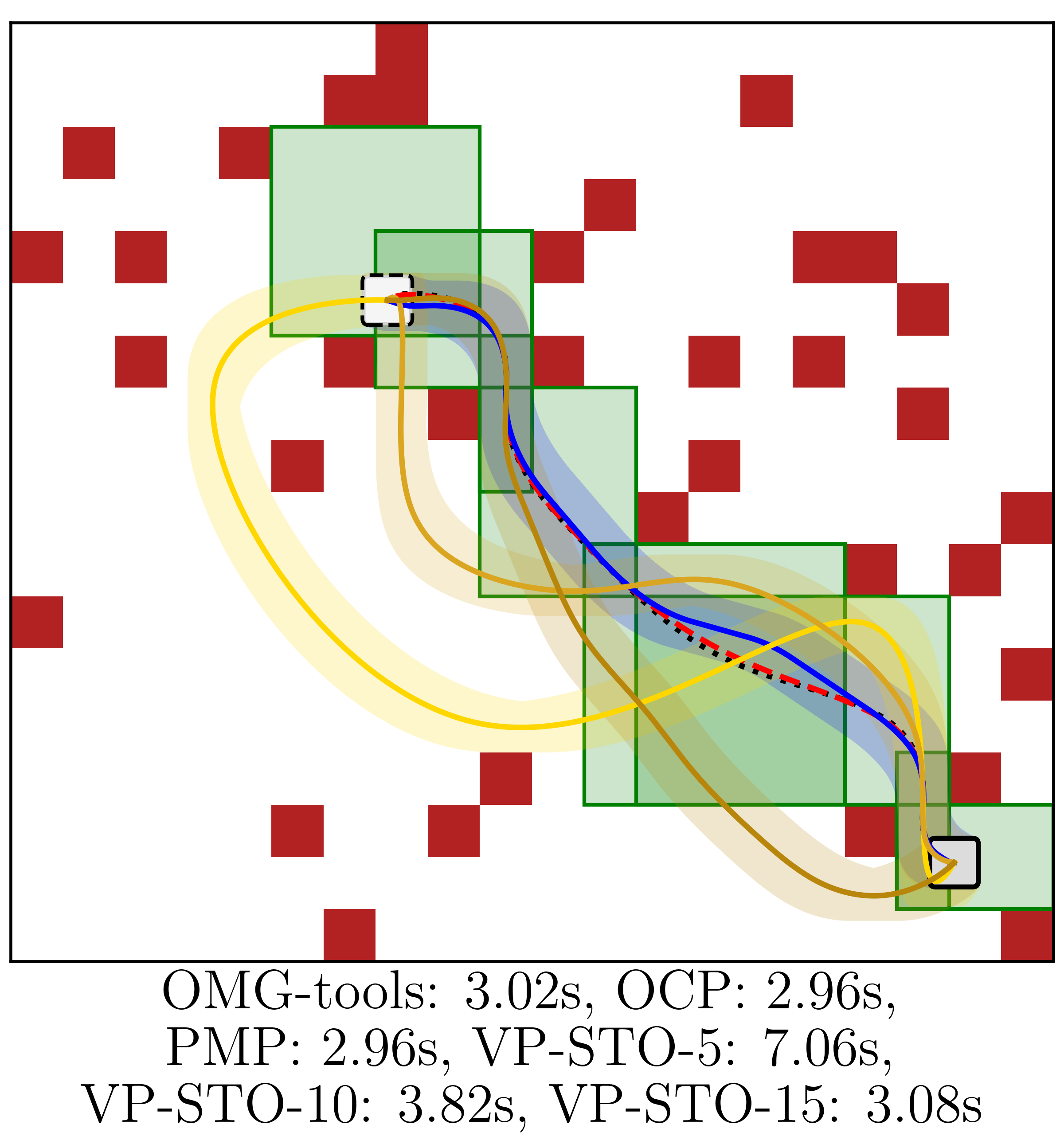}
        \subcaption{}
        \label{fig:unstructured-envs-c}
    \end{subfigure}
    \begin{subfigure}{0.45\linewidth}
        \includegraphics[width=\linewidth]{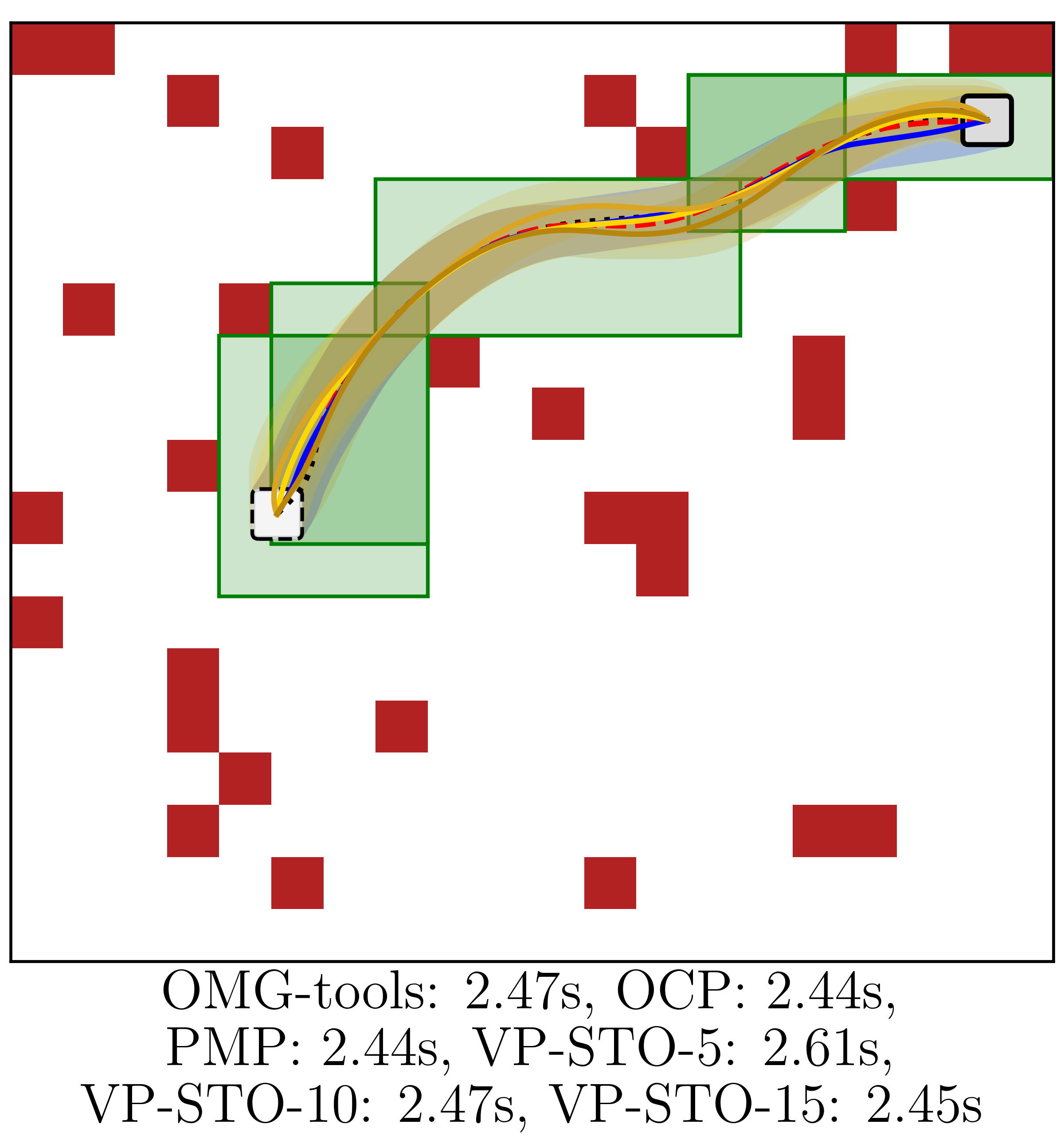}
        \subcaption{}
        \label{fig:unstructured-envs-d}
    \end{subfigure}
    \begin{subfigure}{\linewidth}
        \includegraphics[width=\linewidth]{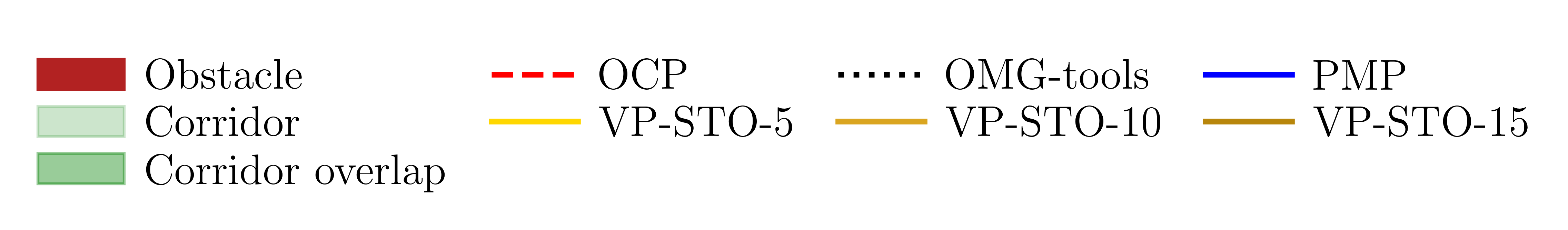}
    \end{subfigure}
    \caption{Some example trajectories in unstructured environments. The vehicle's footprint at $\bm{p}_0$ (solid) and $\bm{p}_n$ (dashed) is shown. Transparantly shaded regions show traces of the vehicle footprint}
    \label{fig:unstructured-envs}
\end{figure}

To evaluate the optimality of the produced trajectories, we consider the computed moving time $t_{\mathrm{move}}$. Table \ref{tab:my_label} shows the moving time, averaged over all trajectories as well median and standard deviation of the relative error of the moving time where OCP is considered the baseline. Figure \ref{fig:suboptimality} shows boxplots and individual samples of this relative error. PMP is less than 0.5\% slower than OCP on average, which is fastest of all methods on average in both types of environments. OMG-tools is around 6\% slower than OCP which is explained by its use of smooth splines to represent the controls. In Figure \ref{fig:suboptimality}, it is clear that OMG-tools is always slightly slower than OCP. PMP mostly converges to the same $t_{\mathrm{move}}$ as OCP, but can produce suboptimal results in some cases.

VP-STO, which is not constrained to a corridor sequence, is 16\% slower on average than OCP for $N_\mathrm{via} = 5$. For $N_\mathrm{via}$ larger than 10, the median relative error on the moving time ($\varepsilon_\mathrm{move}$) is close to zero, indicating that a higher number of via points produces better trajectories. This is also seen in Figure \ref{fig:unstructured-envs-d} where all methods converge to roughly the same trajectory. Still, VP-STO-5 computes a slower trajectory which is caused by the fact that the via points are connected with smooth trajectories which do not explicitly optimize time. To achieve time-optimality, a higher number of via-points is required. The large standard deviation on the relative error on $t_\mathrm{move}$ for VP-STO is explained by the fact that this method often finds different trajectories. Figure \ref{fig:suboptimality} demonstrates that even though VP-STO can converge to the global optimum, this is not guaranteed to happen. In Figure \ref{fig:structured-envs-a} for example, VP-STO converges to different (and worse) local optima than the corridor-based approaches. The probability of VP-STO to converge to the global optimum can be increased by increasing the population size and the initial sampling standard deviation, but this comes at the cost of a higher computation time. In addition, a different number of via-points can influence the convergence, as shown in Figure \ref{fig:structured-envs-b}. Figure \ref{fig:structured-envs-c} shows an example where VP-STO converges to a better local optimum than the corridor-based approaches. Figure \ref{fig:structured-envs-d} shows a case where VP-STO-10 fails to find a feasible solution. Instead, it converges to a locally optimal solution where the collision is minimized by moving straight through the wall. This is a result of a poor initial set of samples which does not contain a trajectory around the wall. Note that in this environment, VP-STO-15 is also unable to find a collision-free trajectory.

Similar observations can be made in the unstructured environments, of which some trajectories are shown in Figure \ref{fig:unstructured-envs}. In Figure \ref{fig:unstructured-envs-a}, the corridor sequence selects a trajectory that involves turns in narrow passages where VP-STO finds a trajectory that is significantly faster. However, in Figures \ref{fig:unstructured-envs-b} and \ref{fig:unstructured-envs-c}, VP-STO struggles to find a trajectory through narrow passages, which in these cases turns out to be faster. This is explained by the fact that the probability of sampling a feasible trajectory through a narrow passage is low.

\subsubsection{Discussion of Robustness}
When checking the trajectories sampled at 100Hz for collisions, it shows that the OCP solution violates corridor bounds 15 times in the structured environment and 7 times in the unstructured environment. This is caused by the fact that constraints are only applied to discrete points. If the number of grid-points per corridor is increased from 30 to 60, only 2 cases with corridor infeasibilities are found in the structured environment and none in the unstructured environment. However, this change more than doubles $t_{\mathrm{solver}}$ on average. Instead of increasing the number of grid-points, an additional margin could be taken into account to satisfy the corridor bounds which would also affect $t_{\mathrm{move}}$. The other methods all satisfy position constraints continuously.

Apart from OCP, all approaches fail to find a solution in some cases. For the corridor-based approaches, specifically PMP and OMG-tools, this happens more often in the unstructured environment. In such cases, the trajectory through subsections of $\bm{C}$ could be computed and stitched together. The trajectory through a single corridor can always be found analytically if the braking distance in each dimension (given by $\left(v^{(.)}\right)^2/\left(2a_{\mathrm{max}}\right)$) does not exceed the available margin to the corridor bounds. If the planner fails and the initial velocity is too high such that the corridor bounds can no longer be avoided, an emergency behaviour is required to make the vehicle come to stop. However, since the proposed method computes full trajectories with zero terminal velocity, vehicles can only encounter such a situation without a previously computed safe trajectory available if the corridor sequence was suddenly changed. For VP-STO, solver failures refer to a solution in which the vehicle moves through obstacles, such as in Figure \ref{fig:structured-envs-d}. These failures occur more often in the structured environment in which large obstacles are present and in which the trajectory must pass through narrow passages. In the unstructured environments considered here, VP-STO tends to converge to a trajectory avoiding narrow passages, as in Figure \ref{fig:unstructured-envs-b}.

\section{Experimental Validation} \label{Sec:experimental-validation}
\begin{figure}
    \centering
    \begin{subfigure}{0.35\linewidth}
        \includegraphics[width=\linewidth]{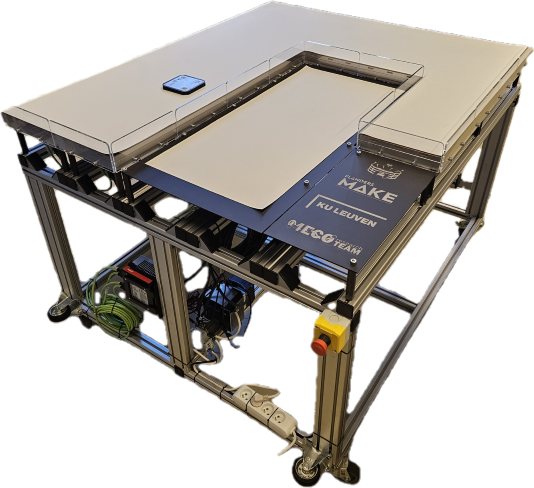}
        \subcaption{XPlanar hardware}
    \end{subfigure}
    \begin{subfigure}{0.6\linewidth}
        \includegraphics[width=\linewidth]{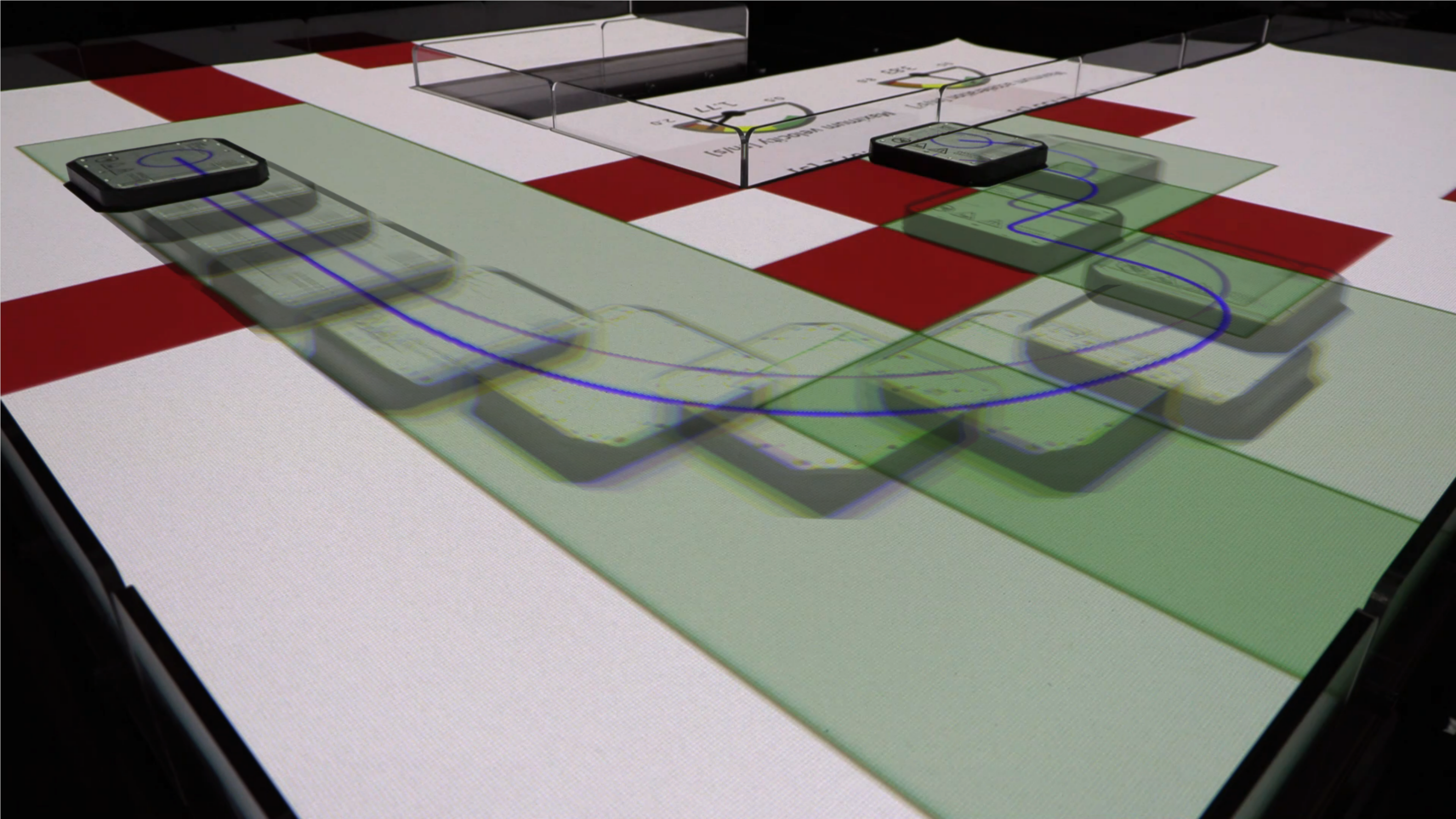}
        \subcaption{Mover executing a trajectory}
    \end{subfigure}
    \caption{Hardware used for experimental validation. Using a projector attached to the ceiling, obstacles, corridors and the trajectory can be visualized live on the hardware}
    \label{fig:hardware}
\end{figure}

We use the XPlanar mover system by Beckhoff Automation as a platform for experimental validation. There are 20 APS4322 XPlanar tiles each measuring 240mm $\times$ 240mm and one APM4220 mover measuring 113mm $\times$ 113mm. The hardware is shown in Figure \ref{fig:hardware}. A video showing a demonstrator on the real system can be found at {\small \url{https://www.youtube.com/watch?v=zhZ2Ko5VxUk}}

The trajectories are sampled at 100Hz and a \texttt{C++} script writes these samples into a FIFO buffer in the PLC using ADS communication. The PLC passes the samples to the \texttt{ExternalSetPointGenerator()} functionality, provided by Beckhoff. The samples contain position, velocity and acceleration setpoints. An internal low-level tracking controller tracks the reference.

\begin{figure}
    \centering
    \includegraphics[width=0.9\linewidth]{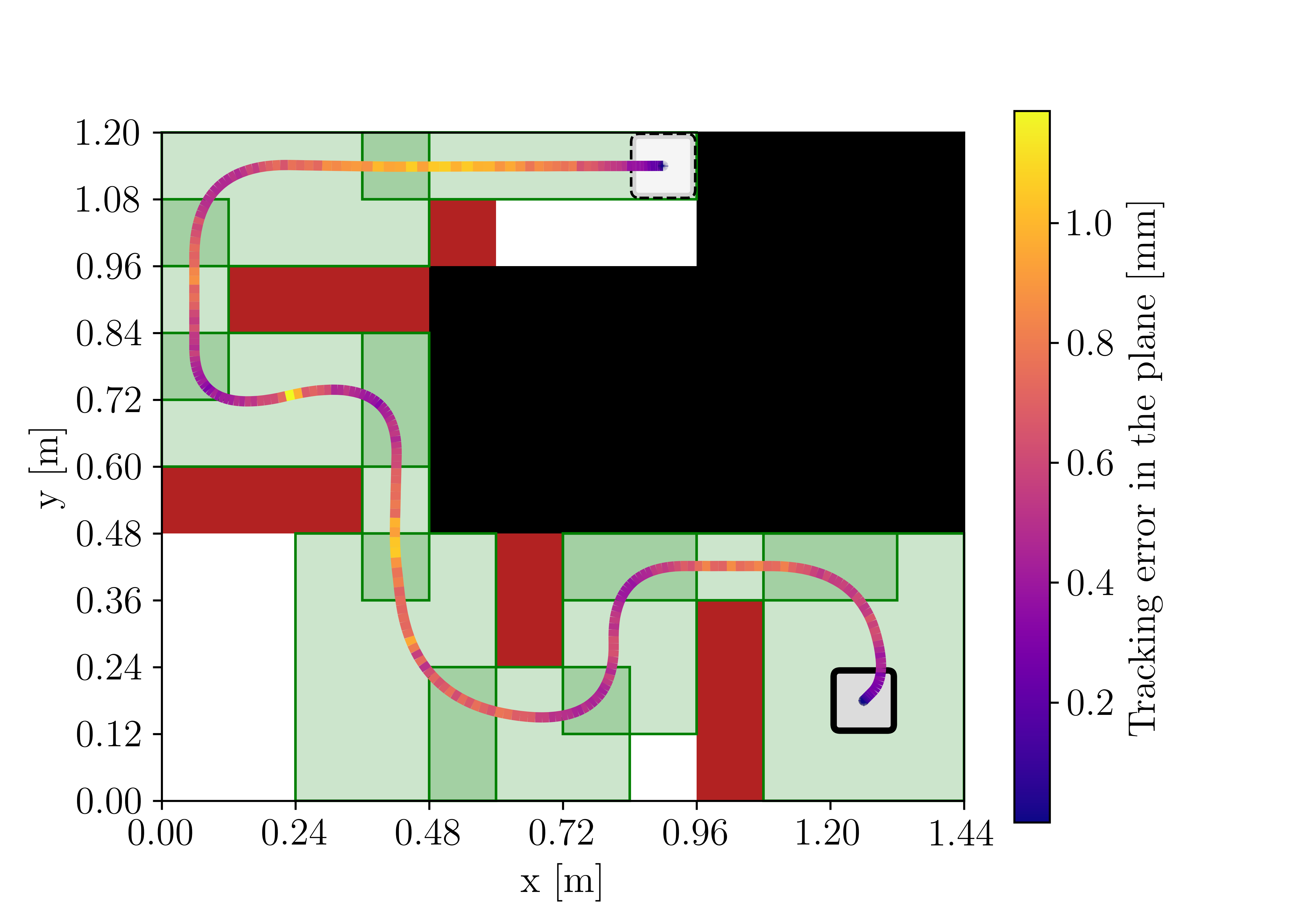}
    \caption{Measured trajectory executed by a real-life mover tracking a reference computed using the proposed method. Red squares represent obstacles. Green rectangles represent the corridor sequence where overlapping areas are shown in a darker shade. The color of the trajectory indicates the distance between the measured position and the provided position setpoint. Note that this error is expressed in mm. The maximum tracking error in the plane is 1.19mm and the average tracking error in the plane is 0.57mm}
    \label{fig:xy-position-tracking-error}
\end{figure}

Figure \ref{fig:xy-position-tracking-error} shows a recorded traveled trajectory of the mover, computed using the proposed method (PMP). For this trajectory, the mover hovered at a height of 3mm, $a_\mathrm{max}$ was set to 6 m/s$^2$ and $v_{\mathrm{max}}$ was set to 2~m/s (which is the maximum velocity according to the technical data). The velocity limit is reached in the $x$-dimension and the maximum value for $|v_y|$ was 1.8m/s. Similar results were found for different parameter values although much higher acceleration limits push the performance of the tracking controller causing the mover to deviate more from the provided setpoints. This can cause the mover going into an error state while stopping abruptly. We consider the obstacles locations to be known and the proposed method is given a map of the environment including these obstacles as mentioned in Section \ref{sec:problem-description}.
The tracking error, meaning the distance between the measured mover position and the desired position setpoint, is indicated by the color of the trajectory. 
The tracking error is always less than 1.2~mm and is higher when moving at higher velocity. This result demonstrates that the vehicle is able to accurately track the trajectories computed by the proposed method and showcases the excellent tracking capabilities of the hardware. This also means that there might be no need to replan trajectories at a fixed rate but rely on the internal reference tracking controller to reach the destination.

\section{Conclusions and Future Work} \label{sec:conclusions}
In this work, we proposed a novel method for computing near time-optimal trajectories for a holonomic vehicle in complex but structured environments. First, we construct a sequence of rectangular corridors to model the free-space. Based on the geometry of these corridors, we use heuristics to determine the locations of waypoints through which we enforce the vehicle to pass as well as the sign of the acceleration at those waypoints. Using this information, we select parametric motion primitives before solving an optimization problem to optimize the remaining degrees of freedom. 

By performing extensive numerical experiments, we show that we can solve the resulting optimization problem significantly faster than state-of-the-art methods, without significantly sacrificing optimality within the corridor sequence. In some cases, however, the corridor sequence is suboptimal, which is made obvious by benchmarking against VP-STO. We plan to address this issue by exploring paths that belong to different homotopy classes and solving different trajectories in parallel. Still, sampling-based approaches only achieve global optimality asymptotically and especially in the case of narrow passages or large obstacles, VP-STO can converge to a (potentially poor) local optimum. This is also influenced by the number of via-points considered, which is a tuning parameter that is difficult to tune. In the proposed method, the number of primitives is automatically fixed by the number of corridors constructed. Future work plans to improve the corridor construction, potentially generating multiple sequences that can be solved in parallel.

We also validated our approach on real, industrial hardware using the XPlanar mover system by Beckhoff Automation. These results indicate the proposed trajectories can be tracked with high accuracy using the low-level controller. 

The logical next step is to use the proposed planner in a multi-agent setting, which will require online replanning. This could for example be done by having vehicles claim areas in the environment and have them plan through these areas, which would fit in a framework related to the one proposed by \citet{devos2024}. Note that because the planner computes time-optimal trajectories, recursive feasibility should be explicitly taken into account by the framework to make sure vehicles that are already moving are not demanded to suddenly evade new obstacles.

The proposed method can be extended relatively easily to multiple dimensions including rotation, since the heuristics and the resulting optimization problem are largely decoupled for every dimension. Additionally, the shape of the motion primitive can be modified to allow for motions with limited jerk. These additions can extend the proposed method to OMRs. 

Finally, it would be interesting to evaluate the performance of the proposed methodology for more general corridor sequences, such as corridors that are slightly tilted and might be able to better represent free-space in some applications. The performance of the heuristics still needs to be evaluated in that case.

\backmatter
\section{Declarations}
\bmhead{Ethics approval and consent to participate} Not applicable.
\bmhead{Consent for publication} All authors declare their consent to publish.
\bmhead{Funding} This work was funded by Flanders Make through the SBO project ARENA (Agile \& REliable NAvigation).
\bmhead{Author's contributions} Louis Callens wrote the main manuscript text and carried out most of the research. Bastiaan Vandewal and Ibrahim Ibrahim helped with developing the algorithm and perform experimental validations and provided significant feedback to the manuscript. Jan Swevers and Wilm Decré provided significant research guidance, helped sketch outlines of the manuscript content and reviewed the manuscript.
\bmhead{Acknowledgments} The authors would like to thank dr. Joris Gillis for his help in using FATROP through CasADi and general \texttt{C++} implementation aspects. Additionally, the authors thank dr. Alejandro Astudillo Vigoya for his feedback on the manuscript.

\bibliography{my-bib}

\end{document}